\documentclass[12pt]{amsart}
\usepackage{amsmath,amscd,latexsym,verbatim,amssymb}
\usepackage[english]{babel}     
\usepackage{times} 

 \newcommand{\resp}{{\it resp.} }
\newcommand{\cf}{{\it cf.} }
\newcommand{\ie}{{\it i.e.} }
\newcommand{\eg}{{\it e.g.} }
\newcommand{\prf}{\noindent {\bf Proof.} }
\renewcommand{\qed}{\hfill$\Box$\medskip}

\newcommand{\sM}{\mathcal{M}}

\newcommand{\sT}{\mathcal{T}}

\newcommand{\sU}{\mathcal{U}}
\newcommand{\sV}{\mathcal{V}}
\newcommand{\sW}{\mathcal{W}}
\newcommand{\sX}{\mathcal{X}}
\newcommand{\sY}{\mathcal{Y}}
\newcommand{\sZ}{\mathcal{Z}}

\newcommand{\C}{\mathbf{C}}
\newcommand{\Q}{\mathbf{Q}}
\newcommand{\Z}{\mathbf{Z}}

\newcommand{\R}{\mathbf{R}}
 \newcommand{\inj}{\hookrightarrow}

\newcommand{\surj}{\rightarrow\!\!\!\!\!\rightarrow}

\renewcommand{\Im}{\operatorname{Im}}

\renewcommand{\epsilon}{\varepsilon}
\font\sm=cmr10 at9pt

\swapnumbers

\newtheorem{thm}{Theorem}[subsection]
\newtheorem{lemme}[thm]{Lemma}
\newtheorem{slemme}[thm]{Sublemma}
\newtheorem{prop}[thm]{Proposition}

\newtheorem{cor}[thm]{Corollary}

\newtheorem{sorite}[thm]{Sorite}
\theoremstyle{definition}
\newtheorem{defn}[thm]{Definition}
\newtheorem{notn}[thm]{Notation}

\numberwithin{equation}{section}
  
 \renewcommand{\qed}{\hfill $\square$\medskip}

\font\big=cmr10 at13pt
\font\sm=cmr10 at11pt
 
\begin{document}  \centerline
{\big\bf ON A GEOMETRIC DESCRIPTION OF $\;Gal({\bar {\Q}_p}/{\Q_p})$,  }

\medskip 
 \centerline
{\big\bf AND A $\,p$-ADIC AVATAR OF $\;\widehat{GT}$}
\title{}
 \author{Yves
Andr\'e}
\address{Institut de Math\'ematiques de
Jussieu\\175--179 rue du
Chevaleret\\ \break 75013
Paris\\France.}
\email{andre@math.jussieu.fr}
\maketitle
 \bigskip\bigskip\bigskip
 \bigskip\bigskip\bigskip

 \tableofcontents

{\sm ABSTRACT. We develop a $p$-adic version of the so-called Grothendieck-Teichm\"uller theory (which studies
$Gal(\bar\Q/\Q)$ by means of its action on profinite braid groups or mapping class groups). For every place
$v$ of $\bar\Q$, we give some geometrico-combinato\-rial descriptions of the local Galois group
$Gal(\bar\Q_v/\Q_v)$ inside $Gal(\bar\Q/\Q)$. We also show that $Gal(\bar\Q_p/\Q_p)$ is the
automorphism group of an appropriate $\pi_1$-functor in $p$-adic geometry.}

\bigskip
{\sm Classification: 11R32, 14H30, 14G22, 14G20, 20F28, 20F36.}

\newpage ${\;}$  \bigskip 
 \section{Introduction} 

\bigskip
\subsection{} Profinite groups which are the absolute Galois group $G_k= Gal(\bar k/k)$ of some field $k$ have
been the object of extensive study, especially in the case where $k$ is a number field. In that case, $G_k$
has more structure: it comes equipped with a constellation of closed subgroups $G_{k_v}=Gal(\bar {k_v}/{k_v})$
attached to the places
$v$ of $\bar k$, and fitting together in an arithmetically relevant way.  

The problem then arises to {\it describe these local Galois groups $G_{k_v}$  
in $G_k$}. 

\subsection{} A purely group-theoretic approach to this problem has been found by Artin (for archimedean $v$)
and Neukirch (for non-archimedean $v$): namely, the subgroups $G_{k_v}\subset G_k$ for $v_{\mid k}$ real are
exactly the subgroups of order two; the subgroups $G_{k_v}$ for $v$ non-archimedean are exactly the closed
subgroups which have, abstractly, the algebro-topological structure of the absolute Galois group of a local
field\footnote{which is ``known", \cf \cite[\S\,VII.5]{nsw}.}, and are maximal for this property, \cf
\cite[\S\,XII.1]{nsw}. This is however more a characterization than a description of the local Galois groups
- rather, of the set of local Galois groups $G_{k_v}$ attached to a fixed place $v_{\mid k}$ of $k$.
   
\subsection{} In this paper, we examine the problem from a completely different viewpoint, aiming at a
{\it geometric} solution (with combinatorial flavour). Our approach is inspired by Grothendieck's
leitmotiv of studying
$G_k$ via its outer action on the ``geometric" algebraic fundamental group  $\pi_1^{alg}(X_{\bar k})$ of
smooth geometrically connected algebraic varieties $X$ defined over
$k\;$. In fact, the problem of describing the local Galois
groups in
$G_k$ has been raised explicitly by Grothendieck in the context of geometric actions on
fundamental groups and his dream of anabelian geometry (\cite{gr}, note $4$)\footnote{we became aware of
these lines on having another look at
\cite{gr} just before completing this work: ``Parmi les points
cruciaux de ce dictionnaire [anab\'elien], je pr\'evois [...] une description des sous-groupes d'inertie de
$\Gamma$ [ $=G_\Q$], par o\`u s'amorce le passage de la caract\'eristique z\'ero \`a la caract\'eristique
$p>0$, et \`a l'anneau absolu $\Z$." }.   
 
\subsection{} In the simplest non-trivial case ($\;k=\Q\;$ and  $\,X={\bf P}^1_\Q\setminus\{0,1,\infty\}\;$),
it turns out that the outer action of $G_\Q$ on the profinite group $\pi_1^{alg}(X_{\bar\Q})$ is faithful
(Belyi). An embedding $\bar\Q \hookrightarrow \C$ being fixed, $\pi_1^{alg}(X_{\bar\Q}) $ is the profinite
completion of the usual `transcendental' fundamental group $\;\pi_1^{top}(X^{^{an}})\;$ of the complex-analytic
manifold
$X^{^{an}}$ (a discrete free group of rank two). 

One can then recover
$\;G_\R\;$ inside $\;G_\Q\;$ as the intersection of
$\;G_\Q\;$ and $\;Out\,\pi_1^{top}(X^{^{an}})\;$ in $\;Out\,\pi_1^{alg}(X_{\bar\Q}) \;$,
(\cf 3.3.1, and 3.3.2 for a more general statement). 
 
\subsection{} We shall give a similar description of the local Galois group $G_{\Q_p}$, embedded into $G_\Q$
via a fixed embedding $\bar\Q \hookrightarrow \C_p=\hat{\bar\Q}_p$. The required ingredient is a
fundamental group for rigid-analytic ``$p$-adic manifolds" playing the role of $\pi_1^{top}$.
In particular, the profinite completion of such a fundamental group should coincide with the algebraic
fundamental group
$\pi_1^{alg}$, in the case of an algebraic $p$-adic manifold.

Such a theory has been set up in \cite{a} (for completely different purposes), and will be outlined
below (\S 4). Let us just say here that its corner-stone is the notion of tempered\footnote{following the
referee's advice, we have changed our previous terminology `temperate' to `tempered'} etale covering, which
generalizes in a minimal way both finite etale coverings and infinite topological coverings.

The corresponding {\it tempered
fundamental groups}
$\pi_1^{temp}$ encapsulate combinatorial information about the bad reduction of all finite etale cove\-rings
of the base. Their topology is a little complicated (for instance, the tempered fundamental group of ${\bf
P}^1_{\C_p}\setminus\{0,1,\infty\}\;$ is complete but not locally compact, hence is neither discrete nor
profinite); but in dimension one, they have a suggestive combinatorial description as inverse limit of a
sequence of fundamental groups of certain finite graphs of groups, which could be considered as {\it
non-archimedean analogues of Grothendieck's `des\-sins d'enfants'}. 

\subsection{} For any smooth geometrically connected algebraic variety $X$ over $\Q_p$, there is a canonical 
outer action of $G_{\Q_p}$ on $\pi_1^{temp}(X^{^{an}})$, where $X^{^{an}}$ denotes the rigid-analytic manifold
attached to $X_{\C_p}$, \cf prop. 5.1.1.

 Coming back to the case of $\,X={\bf P}^1_\Q\setminus\{0,1,\infty\}\;$, one version of our main result is that
{\it one can recover $G_{\Q_p}$ inside
$\;G_\Q\;$ as the intersection of
$\;G_\Q\;$ and $\;Out\,\pi_1^{temp}(X^{^{an}})\;$ in $\;Out\,\pi_1^{alg}(X_{\bar\Q}) \;$} , \cf thm. 7.2.1 (a
more general statement is given in 7.2.3).

  We also study the structure of $Out\, \pi_1^{temp}(X^{^{an}})$ for algebraic curves $X$. We
show that this group is, like $\pi_1^{temp}(X^{^{an}})$ itself, the inverse limit of a sequence of
finitely generated discrete groups ( 6.1.4).

\subsection{} According to the philosophy of the Grothendieck-Teichm\"ul\-ler theory, ${\bf
P}^1_\Q\setminus\{0,1,\infty\}$ should be seen as the special case $n=4$ of the moduli space $\sM_{0,r}$ of
curves of genus $0$ with
$r$ ordered marked points, and one should study as well
the Galois action on the geometric fundamental group of these moduli spaces (and also in higher
genus). The main player here is Drinfeld's Grothendieck-Teichm\"uller group
$\widehat{GT}$, a kind of ``geometric upper bound" for
$G_\Q$ \cite{dr}\cite{i}.

 We introduce a closed subgroup $\widehat{GT}_p$ of
the profinite group $\widehat{GT}$, the {\it local Grothendieck-Teichm\"uller group at $p$}, defined in terms of
the rigid-analytic manifolds attached to the moduli spaces $\sM_{0,r}$ over $\C_p\,$ (8.6.3). This group is a
kind of `` geometric upper bound" for 
$G_{\Q_p}$: we prove that $G_{\Q_p}$ {\it is the intersection of
$\;G_\Q\;$ and $\;\widehat{GT}_p\;$ in} $\;\widehat{GT}\,$ (thm. 8.7.1).

 This sheds some new light on the longstanding problem of how close is $\widehat{GT}$ to $G_\Q$ (\cf 8.7.2).

\subsection{} Another attempt to
describe the absolute Galois group of a number field $k$ 
  as the full automorphism group of some geometric structure consists in looking at sufficiently many
$k$-varieties at a time and at the geometric fundamental group as a functor. 
 Geometric algebraic
fundamental groups give rise to a functor $\underline{\pi}^{alg}_{\bar k}$ from the category $\sV_k$ of smooth
geometrically connected varieties over a number field $k$ to the category $\sT$ of topological groups up to
inner automorphisms. Pop has shown that $G_k= Aut \,\underline{\pi}^{alg}_{\bar k}$ (\cite{p}, unpublished). 

On the other hand, for any $p$-adic place $v$ of $k$, there is a functor $\underline{\pi}^{temp}_{\C_p}:\,\sV_k
\to \sT$ given by the tempered fundamental groups of associated $p$-adic manifolds. {\it Our final result is
that $G_{k_v}= Aut \,\underline{\pi}^{temp}_{\C_p}$} (thm. 9.2.2). 
 
Thus in some sense, the arithmetic of finite extensions of $\Q_p$ is embodied in analytic geometry over $\C_p$.

\newpage
\bigskip \section{Geometric fundamental groups and Galois actions}\label{s1}

\subsection{} Let $X$ be a smooth geometrically connected algebraic variety over over a field $k$, endowed with
a geometric point $x$. Grothendieck's algebraic fundamental group $\pi_1^{alg}(X,x)$ is the profinite group
which classifies all finite etale (pointed) cove\-rings of $(X,x)$:  
\[\{\hbox{finite etale cove\-rings of $X$}\} \sim \{ \hbox{finite }\pi_1^{alg}(X,x)\hbox{-sets}\} .\]
This group depends on $x$ only up to inner automorphism.

 Let $\bar k$ be a separable closure of $k$, and let $G_k=Gal(\bar k/k)$ stand for the absolute Galois group.
The group $\pi_1^{alg}(X_{\bar k},x) $ is sometimes called the geometric fundamental group of $X$ (pointed
at $x$). There is a functorial exact sequence 
\[\{1\}\to \pi_1^{alg}(X_{\bar k},x)\to \pi_1^{alg}(X,x)\to G_k \to \{1\},\]
 which splits canonically if $x$ comes from a $k$-rational point of $X$. 
Whence a Galois action  \[G_k \;\stackrel{\rho}\to  \; Out\,\pi_1^{alg}(X_{\bar k}),\]
which lifts to an action
\[G_k \;\stackrel{ \rho_x}\to \; Aut\,\pi_1^{alg}(X_{\bar k},x)\]
if $x$ comes from a $k$-rational point of $X$ (here $Out= Aut/Inn\,$ denotes as usual the group of outer
automorphisms (in $\rho$, we drop the base point from the notation since it is irrelevant).

\subsection{} We assume henceforth that $ char\, k= 0$. Under this assumption, it is known that the profinite
groups
$\pi_1^{alg}(X_{\bar k},x)$ are finitely generated, which implies that
$ Aut\,\pi_1^{alg}(X_{\bar k},x)$ and $ Out\,\pi_1^{alg}(X_{\bar k})$ are also finitely generated (hence
metrizable
\footnote{recall that a profinite group is metrizable if and only if it has countably many open subgroups, or,
equivalently, if it is a countable inverse limit of finite groups, \cite[\S\,4.1.3]{w}
\cite[\S\,IX.2.8]{bou}, or else, if and only if it is a closed subgroup of a quotient of the profinite
free group on two generators $\hat F_2\,$ \cite[\S\,4.1.6]{w} (a profinite variant of
the fact that any countable discrete group is a subgroup of a quotient of
$F_2$ ); because of these equivalences, one sometimes says `separable' instead of `metrizable' }) profinite
groups. More
precisely, there are only finitely many open normal subgroups of
$\pi_1^{alg}(X_{\bar k},x)$ of index dividing\footnote{another possible choice: replace ``dividing $n$" by
``less or equal to $n$"; our preference is justified \eg by 4.7 below}
$  n$, hence their intersection $U_n$ is a characteristic\footnote{\ie stable under every 
automorphism of the profinite group
$\pi_1^{alg}(X_{\bar k},x)$} open subgroup of
$\pi_1^{alg}(X_{\bar k},x)$ (hence is preserved by $G_k$). It follows that $U_n= \pi_1^{alg}(X_n,x_n)
$ for a well-defined finite Galois etale (pointed) cove\-ring $(X_n, x_n) \to (X_{\bar k},x)$ with $X_n$
geometrically connected; in fact, the tower of characteristic subgroups $U_n$ is the tower of geometric
fundamental groups of a `tower' of (pointed) finite Galois etale cove\-rings 
\[\begin{matrix} \ldots &\to &X_{n'}&\to X_n &\to \ldots &\to X_{\bar k} &\,,\;\;&& n\vert n'.\\ &&
x_{n'}&\mapsto x_n &\mapsto \ldots &\mapsto x &\, \;\;&&  \end{matrix}\] We denote by $\Gamma_n =
\pi_1^{alg}(X_{\bar k},x)/U_n$ the Galois group of the cove\-ring $X_n/X_{\bar k}$. Notice that the above
`tower' is canonical and functorial in $X$, and the same is true for the corresponding `tower' of finite groups 
\[ \ldots\to \Gamma_{n'}\to \Gamma_n\to \ldots \to \{1\} \,,\;\; n\vert n'.\]
(Moreover, its formation is compatible with extension of algebraically closed fields $\bar k\inj \bar k'$).
 One then has
(\cf
\cite[\S\,3.4. ex.6]{w} for the case of $Aut$)
\[ \pi_1^{alg}(X_{\bar k},x)=  \varprojlim \Gamma_n,\]
\[Aut\,\pi_1^{alg}(X_{\bar k},x)=  \varprojlim Aut\,\Gamma_n,\] 
\[Out\,\pi_1^{alg}(X_{\bar k})=  \varprojlim Out\,\Gamma_n.\] 
The topology of the group of (outer) automorphisms defined by the inverse limit is called the {\it
topology of congruence subgroups}; thus, the group of
(outer) automorphisms of
$\pi_1^{alg}(X_{\bar k},x)$ is closed for the (profinite) topology of congruence subgroups.

\medskip\noindent 2.2.1. {\it Observation.}\footnote{I am endebted to one of the referees for this observation}
For any normal closed subgroup $U$ of $\pi_1^{alg}(X_{\bar k})$ contained in $U_n$, the quotient $U_n/U$ is a
characteristic subgroup of $\pi_1^{alg}(X_{\bar k})/U$. 

  Indeed, by definition, the homomorphisms $\pi_1^{alg}(X_{\bar k})\to \pi_1^{alg}(X_{\bar k})/U \to
\Gamma_n$ induce bijections between the sets of finite quotients of order dividing $n$ of each of these
groups. 

In particular, the kernel of $\Gamma_{n'}\to
\Gamma_n$ is characteristic, for any multiple $n'$ of $n$. 
 
\subsection{} In this paper, `curve' is an abbreviation for `smooth geometrically connected algebraic variety of
dimension one'. A curve $X$ over $k$ is called {\it hyperbolic} or {\it anabelian} if its geometric fundamental
group is non-abelian. If $X$ is affine, this just means that $X_{\bar k}$ is not isomorphic to the projective
line minus one or two points.  

The following result will be of constant use throughout the present paper:

\begin{thm}[Belyi, Matsumoto\cite{m}] Assume that $k$ is a number field. Then the Galois action 
\[G_k \;\stackrel{\rho}\to  \; Out \,\pi_1^{alg}(X_{\bar k}) \] is faithful for any hyperbolic affine curve
$X$. \qed
\end{thm}

This relies in turn on Belyi's theorem (used several times in the sequel) according to which any curve defined
over $k$ admits a rational function which is ramified only above $0, 1,\infty$.
 
\bigskip \section{ $G_\Q, G_\R$ and geometric fundamental groups }\label{s2}

\subsection{} Now $k$ is a number field, and we fix an embedding $\iota:\,\bar k\inj \C$. For any smooth
geometrically connected algebraic $k$-variety $X$ with a geome\-tric point $x$, one has canonical isomorphisms
\[\pi_1^{alg}(X_{\bar k},x) = \pi_1^{alg}(X_{\C},x)= \widehat{\pi_1^{top}(X_\C^{^{an}},x)},\] 
the latter group being the profinite completion of the usual topological fundamental group
${\pi_1^{top}(X^{^{an}},x)}$ which classifies topological (= \'etale) cove\-rings\footnote{not necessarily finite}
of the complex-analytic manifold attached to $X_\C$.

If $\iota(k)\subset \R$, there is a Galois action  \[G_\R=\Z/2\Z
\;\stackrel{\rho_\infty}\longrightarrow  \; Out\,{\pi_1^{top}(X^{^{an}} )},\] which lifts to an action
\[G_\R \;\stackrel{ \rho_{\infty,x}}\longrightarrow \; Aut\,{\pi_1^{top}(X^{^{an}},x )}\]
if $x$ comes from a real point of $X$. This action is compatible with the (outer) action of $G_k$ on the
profinite completion  $\pi_1^{alg}(X_{\bar k},x)$.

\subsection{} Assume that $X$ is an {\it affine curve}. Then ${\pi_1^{top}(X^{^{an}},x)}$ is a free
group of finite rank, which implies that it is residually finite, \ie embeds into its profinite completion. It
follows that the following natural homomorphisms are injective:
\[   {\pi_1^{top}(X^{^{an}},x )}\inj  \pi_1^{alg}(X_{\bar k},x),\]
\[  Aut\,{\pi_1^{top}(X^{^{an}},x )}\inj Aut\,\pi_1^{alg}(X_{\bar k},x),\]
\[ Out\,{\pi_1^{top}(X^{^{an}}  )}\inj Out\,\pi_1^{alg}(X_{\bar k}). \]
 For the injectivity of the third homomorphism, one uses in addition the
following 

\begin{lemme} A free group $F$ of finite rank $>1$ is its own normalizer in its profinite completion $\hat F$.
\end{lemme}

For lack of reference\footnote{after submission of this paper, the reference \cite[thm. 2]{gro} was pointed
out to the author by A. Tamagawa and one of the referees}, we indicate a proof. Let
$a$ belong to the normalizer of $F$ in $\hat F$, and let $x_1, \ldots, x_r$ be a basis of $F$. A classical
result of Stebe says that
$F$  is conjugacy-separated: if $x,y\in F$ are
conjugate in every finite quotient of $F$, they are conjugate in $F\,$ \cite[prop. 4.9]{lys}. In particular,
there are elements $a_i\in F\,$ such that $ax_ia^{-1}=a_ix_ia_i^{-1}$. Hence $a_i^{-1}a$ belongs to the
centralizer of $x_i$ in $\hat F$, which is $x_i^{\hat \Z}$. Therefore, there are elements $z_i\in \hat \Z$ such
that $a=a_ix_i^{z_i}$. For $i=1, 2$, this implies $x_1^{z_1}x_2^{-z_2}\in F$, whence $z_1, z_2\in
\Z$, and $a\in F$.\qed
 
\noindent We denote by $\overline{ Out}\,{\pi_1^{top}(X^{^{an}}  )}$ the completion of ${
Out}\,{\pi_1^{top}(X^{^{an}}  )}$ with respect to the topology of congruence subgroups, \ie the closure of
$Out\,{\pi_1^{top}(X^{^{an}}  )}$ in
$ Out\,\pi_1^{alg}(X_{\bar k})$: \[\overline{ Out}\,{\pi_1^{top}(X^{^{an}}  )}=  
\varprojlim \,\Im\,[Out\,{\pi_1^{top}(X^{^{an}}  )} \to Out\,\Gamma_n] \] 
(this limit is clearly a quotient of the profinite completion of $Out\,{\pi_1^{top}(X^{^{an}}  )}$ and embeds
into $\varprojlim \, Out\,\Gamma_n =   Out\,\pi_1^{alg}(X_{\bar k})$ ). 

\medskip Since we shall have to compare repeatedly the (outer) automorphism group of a topological group with
the (outer) automorphism group of its profinite completion\footnote{the profinite completion of a topological
group is the inverse limit of its finite quotients (\ie quotients by {\it open} subgroups of finite index)}, we
include here the 

\begin{sorite} $a)$ Let $\phi:G\to H$ be a surjective homomorphism of topological groups, and let $\hat\phi:
\hat G
\to \hat H$ be the induced homomorphism of their profinite completions. 
\\ Let
$\tau$ be an automorphism\footnote{in the category of topological groups, of course} of $G$, and let $\hat
\tau$ be the corresponding automorphism of $\hat G$. If $\hat\tau $ induces an automorphism of $\hat H$ (via
$\hat \phi$), then $\tau$ induces an automorphism of
$H/\ker (H\to \hat H)$.   

$b)$ Let $F$ be an open subgroup of finite index of a topological group $G$. We assume that $G$ embeds
into its profinite completion $\hat G$. Let
$\tau$ be as before an automorphism of $G$. If $\hat\tau $ preserves the image $\bar F$ of $F$ in $\hat G$,
then $\tau$ preserves $F$.   \end{sorite}

\proof $a)$ Replacing
$G$ and
$H$ by
$G/\ker (G\to
\hat G)$ and
$H/\ker (H\to \hat H)$ res\-pectively, one may assume that they embed into their profinite completion. Let  
$Aut(G,\ker \phi)$ (\resp $Aut(\hat G, \ker \hat\phi)$) denote the group of automorphisms of $G$ (\resp $\hat
G$) which preserve the kernel of $\phi$ (\resp $\hat\phi$). The sorite follows from the equality

\centerline{$Aut(G,\ker \phi)= Aut(\hat G, \ker \hat\phi)\cap Aut\,G,\;$}

 \medskip
$b)$ One has $F=\bar F\cap G$, whence

 \centerline{$Aut(G,F)= Aut(\hat G,   \bar F)\cap Aut\,G.\;$}   \qed

\subsection{} For any hyperbolic affine curve $X$ over a number field $k$ (and a fixed complex embedding
$\iota$ of $k$), we have encountered two closed subgroups of
$Out\,\pi_1^{alg}(X_{\bar k})$: 
$\overline{Out}\,{\pi_1^{top}(X^{^{an}}  )}$ and
$G_k$. We shall study their intersection, starting with the case of ${\bf P}^1\setminus\{0,1,\infty\}$.

\begin{thm} If $X= {\bf P}^1_k\setminus\{0,1,\infty\}\;$, then \[\;Out\,{\pi_1^{top}(X^{^{an}}  )}\cap
\, G_k \,= \overline{ Out}\,{\pi_1^{top}(X^{^{an}}  )}\cap
\, G_k  =\bigg\{^{\hbox{$\, G_\R$ \;  if $\;\iota(k)\subset
\R\,$}}_{\hbox{$\{1\} \;$ otherwise.} } \]
   \end{thm}
  
\prf We already know that $G_\R\subset Out\,{\pi_1^{top}(X^{^{an}}  )}\cap
\, G_k$ if $\;\iota(k)\subset
\R\,$. It remains to show that $\;\overline{ Out}\,{\pi_1^{top}(X^{^{an}}  )}\cap
\, G_k \,\subset \, G_\R  $ (from which it follows that the intersection is trivial if $\;\iota(k)\not\subset
\R\,$). 
 We may and shall assume that $k=\Q$. 

Since ${\pi_1^{top}(X^{^{an}}  )}$ is a free group of rank two, 
$Out\,\pi_1^{top}(X^{^{an}})
 $ is faithfully represented on the abelianization 
$({\pi_1^{top}(X^{^{an}}  )})^{ab}\cong \Z^2 $ according to a
classical theorem of Nielsen (\cf
\cite[\S\,I.4.5]{lys}).
 Hence $Out\,\pi_1^{top}(X^{^{an}})
\cong GL_2(\Z)$ and $\;\overline{Out}\,{\pi_1^{top}(X^{^{an}}  )}$ is a quotient of the
profinite completion $\widehat{GL_2(\Z)}$ and commutes with the subgroup $\{\pm id \}\subset GL_2(\Z) 
\cong\;Out\,{\pi_1^{top}(X^{^{an}}  )}$. This subgroup lies in 
$Out\,{\pi_1^{top}(X^{^{an}}  )} \cap \, G_\Q$, and coincides with $G_\R$. It follows that $\;\overline{
Out}\,{\pi_1^{top}(X^{^{an}}  )}$ $\cap \, G_\Q$ commutes with $G_\R$. But it is well-known that $G_\R$ is its
own centralizer in
$G_\Q$ (\cf \cite[\S\,12.1.4]{nsw}). \qed

For more general affine hyperbolic $k$-curves $X$, we have:   

\begin{thm} Assume that some finite etale cove\-ring of $X_{\bar k}$ admits a non-constant rational function
which omits at least three values.
  Assume also that $\;\iota(k)\subset
\R\,$. Then 
\[\;Out\,{\pi_1^{top}(X^{^{an}}  )}\cap
\, G_k \,=     \, G_\R.\] A fortiori, if the geometric
point $x$ is $k$-rational, then the subgroup of
$G_k$ which stabilizes
$\pi_1^{top}(X^{^{an}} ,x )\subset \pi_1^{alg}(X_{\bar k}, x)\;$ is $\;G_\R$. 
\end{thm} 
 \noindent (Note that the assertion about `Aut' is less precise that the statement about `Out'.) 

 \prf\footnote{after submission of this paper, A. Tamagawa (and one of the referees) pointed to me that this
result can also be deduced from \cite[thm. 1.1, rmk. 2.1]{mt}, even without the assumption that some finite
etale cove\-ring of
$X_{\bar k}$ admits a non-constant rational function which omits at least three values: indeed, let
$\Gamma_{g,r}$ be the relevant mapping class group, viewed as a subgroup of index two in the group
$Out^\ast{\pi_1^{top}(X^{^{an}}  )}$ of outer automorphisms preserving the conjugacy class of the local
monodromy at each puncture, and let $\Gamma_{g,r}$. The assertion follows from two facts: a) 
$G_k\setminus \{1\}$ does not intersect (the closure of) the image of $\Gamma_{g,r}$ in
$Out\,{\pi_1^{alg}(X_{\bar k}  )}$ 
\cite[thm. 1.1]{mt}; b) $Out\,{\pi_1^{top}(X^{^{an}}  )}\cap
\, G_k \subset  Out^\ast{\pi_1^{top}(X^{^{an}}  )} $} From the assumption, there is a finite extension
$k'$  of
$k$ in
$\C$, a finite etale cove\-ring
$Y\to X_{k'}$ and a dominant morphism  \[Y\stackrel{\varphi}\to Z={\bf
P}^1_{k'}\setminus\{0,1 ,\infty\} .\] (We may also assume that $X$, $Y$, $Z$ are endowed with compatible
$k'$-rational geometric points $x, y,z$ respectively). The result
will follow from 3.3.1 by the following d\'evissage:

\begin{lemme} The statement \[(\ast)_{X,k}\;\;\;\;\; { Out}\,{\pi_1^{top}(X^{^{an}}  )}\cap \, G_k \,\subset \,
G_\R  \] satisfies
\\$ a)$  $(\ast)_{X,k'}  \Rightarrow (\ast)_{X,k}\; $ if $k'$ is a finite extension of $k$ in $\C$,   
\\$ b)$ $(\ast)_{Z,k} \Rightarrow (\ast)_{Y,k} \;$ if there is a non-constant morphism $f:\,Y\to Z$,
\\$ c)$  $(\ast)_{Y,k}  \Rightarrow (\ast)_{X,k}\; $ if $Y$ is a finite \'etale cove\-ring of $X$,

 provided 
$\iota(k)\subset \R$. 
 \end{lemme} 
 \prf
  $ a)$  $Out \,{\pi_1^{top}(X^{^{an}}  )}\cap \, G_k $ is finite by $(\ast)_{X,k'}$, and contains $G_\R$
because $\iota(k)\subset \R)$. We conclude by the following 

\begin{slemme} (Artin, \cf
\cite[\S\,12.1.7]{nsw}) If $H$ is a closed subgroup of $G_k$ containing $G_\R$ as a subgroup of finite index,
then $H=G_\R$.\qed
\end{slemme}
  
 For $b)$, it is more convenient to deal with $Aut$ instead of $Out\,$; we pass from the latter to the
former using the following consequence of 2.3.1 (replacing the assumption that $\iota(k)\subset \R$  by
the assumption that
$x$ is $k$-rational, as we may using $a)$ and its trivial converse):

\begin{slemme} The natural homomorphism \[H_X :=(Inn\,\pi_1^{alg}(X_{\bar k}, x). Aut\,
\pi_1^{top}(X^{^{an}},x))
\cap G_k \to Out\, \pi_1^{top}(X^{^{an}}) \cap
G_k\] is an isomorphism (the first intersection is taken in $Aut \,\pi_1^{alg}(X_{\bar
k}, x)\,$, the second one in $\;Out \,\pi_1^{alg}(X_{\bar
k})\,$). \qed \end{slemme}
 
  $b)$ The dominant morphism $f$ induces a homomorphism 
\[\pi_1^{top}(Y^{^{an}},y) \stackrel{f_\ast}\to \pi_1^{top}(Z^{^{an}},z) \] whose image is of finite index. We
set, for short,
\[ G=\pi_1^{top}(Z^{^{an}},z),\;\; F= f_\ast \pi_1^{top}(Y^{^{an}},y),\] and denote by $Aut(\hat G, F)$ the
group of automorphisms of $\hat G$ preserving $F$.

\noindent This homomorphism 
$f_\ast$ induces in turn a homomorphism
\[H_Y\to (Inn\,\hat G.  Aut(\hat G,F))
\cap G_k  \] 
which is injective (because of the intersection with $G_k$ in both terms). By 3.3.4, it thus suffices to show
that 
$(Inn\,\hat G.  Aut(\hat G,F))
\cap G_k$ is contained in  $H_Z= (Inn\,\hat G.  Aut\,G)
\cap G_k$. We shall prove the stronger statement \[Aut(\hat G,F)
\subset Aut\,G.\]
 To this aim, it will be convenient to replace $F$ by a subgroup $F'$ which is
still stable under $ Aut(\hat G, F)$ but has the advantage of being {\it normal} (and of
finite index) in
$ G$. We take
 \[\displaystyle F'= \bigcap_{g,\gamma}\, Inn(\gamma(g))
(F)=\bigcap_{\gamma}\bigcap_{g} \gamma(Inn(g)F) 
 =\bigcap_{\gamma}\gamma(\bigcap_{g} \,Inn(g)F),\]
where the intersection runs over $g\in G, \gamma \in Aut(\hat G,F)$. This is clearly a normal subgroup of
$G$, stable under $ Aut(\hat G,F)$. To see that it is of
finite index, note that each $\gamma(\bigcap_{g}
\,Inn(g)F)$ is a subgroup of the same (finite) index in $
G$, and there are only finitely many such subgroups since 
$\hat G$ is of finitely type (\cf \cite[\S\,4.1.2]{w}). 

It thus suffices to show that $Aut(\hat G, F')\subset  Aut\,G$, or equivalently, that
the subgroup $F''$ of $\hat G$ generated by elements of the form $\gamma(g), \; g\in
G, \,\gamma \in Aut(\hat G, F')$ is nothing but $G$. Note that
$F'$ is normal in $F''$.
 There are exact sequences $\{1\}\to   F' \to G\to \Gamma\to \{1\}, \;\;\{1\}\to 
\widehat{F'}\to \hat G\to \Gamma\to \{1\} $ with $\Gamma
$ finite, whence an exact sequence
\[\{1\}\to F''\cap  \widehat{F'}\to F''\to \Gamma\to \{1\} .\]
But $F'$ is normal in the subgroup $F''\cap  \widehat{F'}$ of its profinite completion $\widehat{F'}$.
Since
$F'$ is free of rank $>1$, we conclude from 3.2.1 that $F''\cap \widehat{F'} 
 =F'\subset G$, hence $G=F''$ since $G$ surjects onto $\widehat{G}/\widehat{F'}=\Gamma$.

\medskip $c)$ By the previous step, we may replace $Y/X$ by a bigger finite etale covering, and thus assume that
$\pi_1^{alg}(Y_{\bar k},y)$ is a characteristic subgroup of $\pi_1^{alg}(X_{\bar k},x)$ with finite quotient
group $\Gamma$ (the same property then holds for the topological fundamental groups). The canonical homomorphism
\[Aut\,\pi_1^{alg}(X_{\bar k},x)\to Aut\,\pi_1^{alg}(Y_{\bar k},y)\] has finite kernel, and induces a
homomorphism
\[Out\,\pi_1^{alg}(X_{\bar k} )\to (Out\,\pi_1^{alg}(Y_{\bar k} ))/\Gamma\] which also has finite kernel.
Hence $Out\,\pi_1^{top}(X_{\bar k} )\cap G_k$ is finite if $Out\,\pi_1^{top}(Y_{\bar k} )\cap
G_k$ is, and one concludes by 3.3.4.   
\qed

  \medskip\noindent {\it Remarks.} 1) Some heuristic arguments of Baire type on moduli spaces for curves with
marked points seem to indicate that not every affine hyperbolic curve satisfies the assumption of 3.3.2 (it
would be nice to have a complete proof)\footnote{I am grateful to L. Ramero and P.
Colmez for interesting discussions on this topic}.  
 
 2) The assumption of 3.3.2 is fulfilled for $X= $ an elliptic curve $ E$ minus one point. Indeed, by
translation, one may assume that this point is the origin $O$. Then $E$ minus $E[2]$ (the $2$-torsion
points) is an abelian etale covering of $E\setminus \{O\}$, and the quotient of $E\setminus  E[2]$ by the
involution $P\mapsto -P$ is isomorphic to  ${\bf P}^1\setminus\{0,1,\infty, \lambda\}$. 

\noindent On the other hand, the argument of 3.3.1 applies directly to $X= E\setminus \{O\}$ since
$\pi_1^{top}(X^{^{an}})\cong F_2$. Note however that whereas $Out\,\pi_1^{top}(X^{^{an}})
 $ is faithfully represented on the abelianization 
$({\pi_1^{top}(X^{^{an}}  )})^{ab}= H_1(X^{^{an}},\Z),$ it does not amount to the same to take the
intersection with the image of $G_k$ in $Out \,\pi_1^{alg}(X_{\bar k})$ or in 
$Out \,({\pi_1^{alg}(X_{\bar k})})^{ab}= GL\, (H_1^{et}(X_{\bar k},\hat\Z))$ (the latter is much bigger in
general). 

 3) The statements of 3.3.1 and 3.3.2 do not change if one replaces $\iota$ by its complex conjugate; hence it
is enough to fix an archimedean place of $\bar k$ (instead of $\iota$).

\bigskip \section{ The tempered fundamental group of a $p$-adic manifold}\label{s3}

\subsection{} We fix a prime number $p$ and denote by $\C_p$ as usual the completion of a fixed algebraic
closure $\bar\Q_p$ of $\Q_p$. Let
$K$ be a complete subfield of 
$\C_p $, and let 
 $\bar K$ be its algebraic closure in $\C_p$. We notice that
$G_K= Gal(\bar K/K)$ is a metrizable profinite group (and even a finitely generated profinite group if $K=\C_p$
or a finite extension of $\Q_p$).

 Among the several approaches to
analytic geometry over
$K$, we have found Berkovich's one
\cite{be} most convenient for a discussion of unramified cove\-rings and fundamental groups; indeed, Berkovich's
spaces are locally ringed spaces in the usual sense, as opposed to Tate's rigid spaces which are
Gro\-then\-dieck's sites (Berkovich's spaces contain `more points' than rigid spaces, but unlike the passage
from classical algebraic varieties to schemes, all Berkovich non-classical points are closed). 

\medskip 
In the sequel, by {\it (analytic) $K$-manifold} $\;$ - or {\it $p$-adic manifold} if $K=\C_p\,$- $\;$  we shall
mean a {\it smooth paracompact strictly $K$-analytic space} in the sense of Berko\-vich \cite{be'}\footnote{in
\cite{a}, we used this term in a less restrictive sense, allowing a non-empty boundary}. Any smooth algebraic
$K$-variety $X$ gives rise to an analytic $K$-manifold $X^{an}$, its analytification.

According to Berkovich \cite{be''}, $K$-manifolds are locally compact, locally arc\-wise connected, and
{\it locally contractible}, hence subject to the usual theory of universal cove\-rings and topological
fundamental groups.

By geometric point or base point of a $K$-manifold, we mean a point defined over some algebraically closed
complete extension of
$K$.  

\subsection{} In the sequel, we shall have to consider etale coverings of possibly infinite degree. Let us
recall a definition \cite{dj}\cite{a} which applies both to the complex and to the non-archimedean situations:
$\sY\to \sX$ is an {\it etale covering} (\resp a {\it topological covering}) if, locally on the base $\sX$, it
is a disjoint union of finite etale covering maps (\resp isomorphisms). 

For complex manifolds, any \'etale cove\-ring gives rise to a topological cove\-ring, and
conversely. 

On the other hand, for non-archimedean $K$-manifolds, topological cove\-rings
give rise to etale cove\-rings \cite[\S \,2.6]{dj}, but not conversely: for instance, the `Kummer cove\-ring'
\[{\bf P}^1_K\setminus\{0, \infty\} \to {\bf P}^1_K\setminus\{0, \infty\},\; z\mapsto z^n,
\,n>1\] is a finite etale cove\-ring, but not a local homeomorphism of analytic $K$-manifolds (in
fact the cardinality of the fiber, drops to $1$ at some non-classical points; in the more traditional
rigid viewpoint, this is because there is no {\it admissible} open cover over which the Kummer cove\-ring
splits).  
 
Besides the topological fundamental group $\pi_1^{top}(\sX,x)$, which classifies as usual the
topological cove\-rings of the pointed connected $K$-manifold $(\sX ,x)$, there is the etale
fundamental group $\pi_1^{et}(\sX,x)$, introduced by de Jong \cite{dj}, which classifies the (possibly infinite)
etale cove\-rings of $\sX$. In fact $\pi_1^{top}(\sX,x)$ is a discrete quotient of the topological group
$\pi_1^{et}(\sX,x)$.

But neither group is a close analogue of the fundamental group of a complex manifold:

$\bullet\; \pi_1^{top}(\sX,x)$ is ``too small": for instance, $ \pi_1^{top}({\bf
P}^1_{\C_p}\setminus\{0, 1,\infty\},x) =\{1\}$. 

$\bullet\; \pi_1^{et}(\sX,x)$ is ``too big": for instance, $ \pi_1^{et}({\bf
P}^1_{\C_p} ,x)$ is non-trivial (in fact, it is a huge non-abelian group \cite[\S \,7]{dj}). 

 \subsection{} In order to remedy this, we have introduced in \cite{a} an intermediate fundamental group
$\pi_1^{temp}(\sX,x)$, the {\it tempered fundamental group}, which classifies all tempered (pointed)
cove\-rings of $(\sX,x)$:  
\[\{\coprod\hbox{tempered cove\-rings of $\sX$}\} \sim \{ \hbox{discrete }\pi_1^{temp}(\sX,x)\hbox{-sets}\} .\]
By definition, an etale cove\-ring $\sY\to \sX$ is said to be {\it tempered} if there is a commutative diagram of
etale cove\-rings
 \[\begin{matrix} &&\sZ&&\\ &\swarrow &&\searrow &\\ \sT&&&& \sY\\& \searrow &&\swarrow& \\&& \sX &&
\end{matrix}\] 
where $\sZ\to \sT$ is a (possibly infinite) topological cove\-ring, and $\sT\to \sX$ is a finite etale cove\-ring. 

In some
sense, this is the `minimal' theory of etale cove\-rings which {\it takes into account both the topological
cove\-rings and the finite etale cove\-rings}.

We refer to \cite[\S\S\,  1, 2]{a} for a discussion of these cove\-rings and the precise definition of
$\pi_1^{temp}(\sX,x)$ as a separated prodiscrete topological group. This topological group depends on the base
point $x$ only up to inner automorphism
\cite[\S\,1.4.4]{a}.   

Here, we shall content ourselves with the following useful `criterion':

\begin{prop}\cite[\S\,2.1.9]{a}  Let $g\,:\sY \to \sX$ be a Galois etale covering with (discrete) group
$G$. If $G$ is torsion-free (\resp virtually torsion-free\footnote{recall that a discrete group is said
to be virtually torsion-free is it admits a torsion-free subgroup of finite index}), then
$g$ is a topological covering (\resp a tempered covering). \qed
 \end{prop} 
 
\medskip\noindent {\it Remarks.} 1) Because of this proposition, infinite covering maps are often easier to
handle than finite ones, when the point is to check whether they are local homeomorphisms; this may justify
the detour to infinite coverings even if one is primarily interested in finite ones.  

2) We mention in passing that there is a theory of tempered
fundamental groups based at tangential base points, in dimension one \cite[\S\,2.2]{a}.  

\subsection{} Let us consider the case when $\sX= X^{^{an}}$ is the analytification of a geometrically
connected smooth algebraic $K$-variety.
 
\begin{prop} \cite[\S\,2.1.7]{a} $\,\pi_1^{temp}(X^{^{an}},x)$ is a
countable inverse limit of discrete finitely generated groups; it particular, it is a polish group. \\ Its
profinite completion is canonically isomorphic to
$\pi_1^{alg}(X,x)$.\qed
\end{prop}
(the proof of the second assertion relies on the Gabber-L\"utkebohmert version of Riemann's existence theorem
\cite{lu}).
 
Let us recall that, according to Bourbaki \cite[\S\,IX.6.1]{bou}, a topological group is {\it polish} if it is
metrizable, complete, and countable at infinity (it might not be locally compact). These form a nice category of
topological groups to work with: let $\{1\}\to N\to G\to H \to \{1\}$ be a sequence of homomorphisms of
topological groups which is exact in the abstract sense. Assume that $N$ is closed in $G$. If $G$ is polish,
then so are
$N$ (and any closed subgroup of $G$) and $G/N$ \cite[\S\,IX.2.8 prop. 12; \S\,IX.3.1 prop.4 ]{bou}), and if $H$
is also polish, then the bijective homomorphism $\,
G/N\to H$ is an isomorphism  
\cite[\S\,IX.5. ex. 28]{bou}. 

\subsection{} In the case $K=\C_p$, $\pi_1^{temp}(X^{^{an}},x)$ can be described as follows. Let us consider
the `tower' of finite Galois etale cove\-rings 
\[ \ldots\to X_{n'}\to X_n\to \ldots \to X \,, n\vert n'\;, \;\;\;Gal(X_n/X)=\Gamma_n,\]
introduced in 2.2. Let us denote by $\tilde X_n$ the universal topological cove\-ring of $X_n^{^{an}}$. Then
$\tilde X_n$ is an etale Galois cove\-ring of $X^{^{an}}$ \cite[2.1.2]{a} with Galois group
$\Delta_n$ sitting in a extension
\[\{1\}\to \pi_1^{top}(X_n^{^{an}}, x_n)\to \Delta_n\to \Gamma_n\to \{1\}, \] where the discrete group
$\pi_1^{top}(X_n^{^{an}},x_n)$ is finitely generated \cite[\S\,1.1.3]{a}. When $n$ increases (by
divisibility), these exact sequences form an inverse system, which is canonical and functorial in $X$. 

One has
\[\pi_1^{temp}(X^{^{an}},x)= \varprojlim \,\Delta_n    \] (surjective transition maps), and there is a canonical
morphism commutative diagram
\[\begin{matrix} \{1\}&\to &\pi_1^{temp}(X_n^{^{an}})&\to &\pi_1^{temp}(X^{^{an}})&\to &\Gamma_n&\to &\{1\}\\
\downarrow &&\downarrow &&\downarrow &&\downarrow &&\downarrow \\ \{1\}&\to &\pi_1^{top}(X_n^{^{an}})&\to
&\Delta_n&\to &\Gamma_n&\to &\{1\} \end{matrix}\]
with surjective vertical maps, \cf \cite[\S\,2.1.5.,
2.1.8]{a}. 

 On the other hand, {\it in dimension one}, $\pi_1^{top}(X_n^{^{an}})$ has a simple combinatorial description
as the fundamental group of the (dual) graph of incidence of the semistable reduction of $X_n$ in characteristic
$p$,
\cf
\cite[\S\,5.3.]{dj}, in particular it is a free group (notably, $X_n^{^{an}}$ is simply connected if it has good
or tree-like reduction). It follows that
$\Delta_n$ is virtually free (\ie admits a free subgroup of finite index), hence residually finite, which
implies that
$\pi_1^{temp}(X^{^{an}},x)$ itself is {\it residually finite}, so that the following natural homomorphism is
injective \cite[\S\,2.1.6]{a}:
\[   {\pi_1^{temp}(X^{^{an}},x )}\inj  \pi_1^{alg}(X_{\bar k},x).\] 

\noindent{\it Remark.} By the theorem of Karass-Pietrovsky-Solitar \cite{kps}, each $\Delta_n$, being virtually free
and finitely generated, is the fundamental group of a finite graph of finite groups. It would be very
interesting to exhibit a canonical geometric construction of such graphs of groups, which could be considered
as $p$-adic analogues of Grothendieck's ``dessins d'enfants". At present, such construction have been proposed
only in the special case where the projective completion of $X_n$ is a so-called Mumford curve, \ie has a
maximally degenerate reduction
\cite{he}\cite{kat}, using the tree of $SL_2$ over a local field\footnote{after submission of this paper and
discussion with F. Kato, it seems that he is now able to construct such $p$-adic dessins d'enfants
in the general case}.

\subsection{} From this description, $\pi_1^{temp}(X^{^{an}})$ can be easily computed when $X$ is a
non-hyperbolic curve. For instance, 

$\bullet\;\; \pi_1^{temp}({\bf
P}^1_{\C_p}\setminus\{0, \infty\}) \cong \hat\Z,$

\noindent and if $X=E_j$ is an elliptic curve with invariant $j$,

$\bullet\;\; \pi_1^{temp}(E_j^{^{an}}) =\pi_1^{alg}(E_j) \cong\hat\Z \times \hat \Z  \;\;$ if $\vert
j\vert_p \leq 1\;$   (and $\Delta_n \cong (\Z/n\Z)^2$),

$\bullet\;\; \pi_1^{temp}(E_j^{^{an}})   \cong \Z \times \hat \Z\;\;$ if $\vert j\vert_p >
1$, \ie if $E_j$ has bad reduction 

\noindent (in this case, one has a commutative diagram 
\[\begin{matrix} \pi_1^{top}(X_n^{^{an}})&\to& \Delta_n&\to &\Gamma_n\\ \downarrow \cong &&\downarrow
\cong &&\downarrow\cong
\\ n\Z&\to& \Z \times (\Z/n\Z) &\to &(\Z/n\Z)^2 &&) .\end{matrix}\]

\medskip
 When $X$ is hyperbolic, there is no such simple explicit description; in general, $
\pi_1^{temp}(X^{^{an}},x)$ is `{\it lacunary}' in the sense that its open subgroups of finite index have many
infinite discrete quotients. In fact, $\,\pi_1^{temp}(X^{^{an}},x)$ encapsulates the combinatorial information
about the reduction of all finite etale cove\-rings of $X$. 

$\bullet\;$ For $X= $ the projective line minus
$n\geq 3$ points,  
$\pi_1^{temp}(X^{^{an}})$ is not locally compact, and depends on the position of the missing points if $n\geq
4$ \cite[\S\, 2.3.12, 4.5.5]{a}). 

\subsection{} The subtle lacunary properties of $\pi_1^{temp}$, on which the sequel depends, really
belong to the ``profinite theory" and are lost if one passes to the maximal prime-to-$p$ quotient or to the
pro-$p$ completion; indeed:   

\begin{prop} If $X$ is a curve with good reduction (if $X$ is affine, we also require that there is no
confluence of the points at infinity by reduction), the natural homomorphism
\[\pi_1^{temp}(X^{^{an}},x)  \to  \pi_1^{alg}(X,x)^{(p')}\times 
\pi_1^{alg}(X,x)^{(p)}\] to the product of the maximal prime-to-$p$ quotient and the maximal pro-$p$ quotient of
$\pi_1^{alg}(X,x)$, is {\rm surjective}. 
\\  
\end{prop} 

\proof Recall that $\Gamma_n$ is the quotient of $\pi_1^{alg}(X)$ by the intersection of its (finitely many)
open subgroups of index dividing $n$. It follows that
$\vert \Gamma_n\vert$ divides a power of $n$, and that there are sujective
homomorphisms \[\Gamma_n^{(p)}\surj \Gamma_{p^{m_1}},\;\; \Gamma_{p^{m_2}}\surj \Gamma_n^{(p)},\] where
$p^{m_1}$ (\resp $p^{m_2}$) is the greatest power of
$p$ dividing
$n$ (\resp $\vert \Gamma_n\vert$). From the equality
$\pi_1^{alg}(X,x)  =\varprojlim 
\,\Gamma_n$, we derive that
\[\pi_1^{alg}(X,x)^{(p')}
=\varprojlim_{(p,n)=1}\, \,\Gamma_n,\;\;\pi_1^{alg}(X,x)^{(p)}
=\varprojlim_m\, \,\Gamma_{p^m} .\]  Since $\pi_1^{temp}(X^{^{an}},x)=\varprojlim  \,\Delta_n$, it
thus suffices to show that if
$n$ is prime to $p$ (\resp is a power of $p$), then  
$\Delta_n=\Gamma_n$, which means that 
$X^{^{an}}_n$ is simply-connected. By construction, $X_n$ is a finite etale Galois
covering with group $\Gamma_n$ of the curve with good reduction $X$, and $\vert\Gamma_n\vert $ is prime to $p$
(\resp
 is a power of $p$). A standard descent argument shows that
$X_n\to X$ is actually defined over a complete discretely valued subfield $K\subset \C_p$ with
algebraically closed residue field. By Grothendieck's specialization theorem \cite{gr'} (\resp Raynaud's
specialization theorem \cite[thm.1']{r}),
$X_n$ has good reduction (\resp tree-like reduction), hence
$X^{^{an}}_n$ is simply-connected.\qed   

\noindent {\it Remark.} It follows from this that $\pi_1^{temp}({{\bf
P}^1_{\C_p}}^{^{an}}  \setminus\{0,1,\infty\})$ maps surjectively onto the free pro-nilpotent group of rank
two (in contrast, it is well-known that $\pi_1^{top}({{\bf
P}^1_{\C}}^{^{an}}  \setminus\{0,1,\infty\})$ maps injectively into the free pro-nilpotent group of rank
two). 

 \bigskip \section{ A fundamental exact sequence}\label{s4}

\subsection{} In this section, we establish an exact sequence relating the `arithmetic' and `geometric'
tempered fundamental groups.

As above, $K$ is a complete subfield of $\C_p$. For any analytic $K$-manifold $\sX$, $\sX_{\C_p}$ has finitely
many connected components
\cite[\S\,2.14]{dj}, and $\sX$ is said to be geometrically connected if there is just one; it amounts to the
same to require that for every finite extension $L/K$, $\sX_L$ is connected.

Let $(\sX,
x)$ be a geometrically connected
 $K$-manifold endowed with a geometric
point.  
Let $L/K$ be a finite Galois extension contained in $\C_p$.
There is a functorial exact sequence of separated prodiscrete
groups \cite[\S\,2.1.8]{a})
\[ \{1\}\to \pi_1^{temp}(\sX_L,x)\to \pi_1^{temp}(\sX ,x)\to Gal(L/K) \to \{1\},\]
 which splits canonically if $x$ comes from a $K$-rational point of $X$.   
When the Galois extension $L/K$ varies, these exact sequences are compatible in an obvious sense, and  
on passing to the limit, they provide another exact sequence of separated prodiscrete
groups 
  \[  \{1\}\to {\varprojlim}_L \pi_1^{temp}(\sX_L ,x)\to \pi_1^{temp}(\sX ,x)\to G_K
\to \{1\}\]
(the surjectivity of $\pi_1^{temp}(\sX ,x)\to G_K$ can be proved exactly as in \cite[\S\,2.12]{dj}). 

\begin{prop} The canonical homomorphism \[\phi:\;\pi_1^{temp}(\sX_{\C_p} ,x) \to
{\varprojlim}_L
\pi_1^{temp}(\sX_L ,x) \] is an isomorphism. Therefore there is a functorial exact sequence  \[  \{1\}\to
\pi_1^{temp}(\sX_{\C_p} ,x)\to \pi_1^{temp}(\sX ,x)\to G_K
\to \{1\},\]  whence a Galois action  \[G_K \;\stackrel{\rho}\to  \;
Out\,\pi_1^{temp}(\sX_{\C_p} ),\] which lifts to an action
\[G_K \;\stackrel{ \rho_x}\to \; Aut\,\pi_1^{temp}(\sX_{\C_p} ,x)\]
if $x$ comes from a $K$-rational point of $X$.
\end{prop}

\proof We shall use the following
  
\begin{lemme} Let $\phi:\,G\to H$ be a (continuous) homomorphism of sepa\-rated prodiscrete groups. Let
$\phi^\ast$ be the induced functor (in the opposite direction) between  
discrete $H$- and
$G$-sets. If $\phi^\ast$ is essentially surjective (\resp an equivalence), then $\phi$ is
injective (\resp an isomorphism). 
\end{lemme}
\cf \cite[\S\,1.4.9]{a}. \qed

In order to prove the proposition, it thus suffices to show that $\phi^\ast$ is an equivalence. It is clearly
fully faithful, and the essential surjectivity amounts to saying that any tempered covering of $\sX_{\C_p}$ is
defined over some finite extension of $K$. This is clear for finite etale coverings, and it suffices to check
that the universal topological covering $\widetilde{\sY_{\C_p}}$ of any finite covering ${\sY_{\C_p}}$ of
$\sX_{\C_p}$ is defined over some finite extension of $K$. 

Let us first record the following easy 

\begin{slemme} Let $\sY$ be a geometrically connected $K$-manifold.
\\ $i)$ For any finite extension $L/K$, $\widetilde{\sY_L}$ is a component of $(\tilde \sY)_L$.
\\ $ii)$ The connected components of $(\tilde \sY)_{\C_p}$ are defined over some finite extension $K'/K$
(\cite[\S\,2.14]{dj}), and then $(\widetilde{\sY_{K'}})_{\C_p}$ is connected.
\\ $iii)$ For any finite extension $L/K'$, $\widetilde{\sY_L}=  (\widetilde{\sY_{K'}})_L$: topological coverings
of $\sY_L$ come from topological coverings of $\sY_{K'}$. \qed
 \end{slemme}  
The fact that the universal topological covering of ${\sY_{\C_p}}$ is defined over some finite extension of $K$
(in fact over the extension $K'/K$ just introduced) follows from the following lemma (applied to
$\sT=\widetilde{\sY_{K'}}$):
 
\begin{lemme} Let $\sT$ be a geometrically connected and simply connected $K$-manifold. Then
$\sT_{\C_p}$ is simply connected (in the sense of coverings). 
\end{lemme}

 \proof Let $h: \sZ\to \sT_{\C_p}$ be a topological covering: every point $t\in \sT_{\C_p}$
(classical or not) has an open neighborhood $V(t)$ such that $h^{-1}V(t)$ is isomorphic to a disjoint sum of
copies of $V(t)$ via $h$. By the argument of \cite[ \,p. 103]{be'}, there is a finite Galois extension $L/K$
contained in $\C_p$ and an open domain $U_L(t)\subset \sT_L$ such that $t\in (U_L(t))_{\C_p} \subset V(t)$.
The restriction of $h$ above $(U_L(t))_{\C_p}$, being a split covering, admits a canonical $L$-structure.
Saturating with respect to $Gal(L/K)$, one finds by descent an open domain $U(t)\subset \sT$ such that $U(t)_L=
\cup_{\sigma \in Gal(L/K)} \,U_L(t)^\sigma$ and such that the restriction of $h$ above $ U (t) _{\C_p}$
admits a canonical $K$-structure. By 5.1.3.$iii)$, this actually defines a topological covering
$\sZ(t)\to U (t) $. These topological coverings glue together (using \cite[\S\,1.3.3.a]{be'}) and define a
topological covering of
$\sT$ which is a $K$-structure for $h$, and which splits since $\sT$ is simply connected (in the sense of
coverings).  
\qed 
 
\noindent {\it Remark.} It does happen that
$(\widetilde{\sY})_{\C_p}$ has several connected components (for a geometrically connected $K$-manifold $\sY$).
The simplest example is given by a twisted Tate elliptic curve $\sY$, \cf
\cite[\S\, 4]{be}. 

\subsection{} Let us assume that $\sX=X^{^{an}} $ for some geometrically connected smooth algebraic $K$-variety
$X$. Then, passing to the profinite completion, one gets a commutative diagram 
\[\begin{matrix}  \{1\}&\to &
\pi_1^{temp}( X_{{\C}_p}^{^{an}},x)&\to &\pi_1^{temp}( X^{^{an}},x)&\to & G_K
&\to &\{1\}\\  & &\downarrow &&\downarrow &&\downarrow && \\ &&
\widehat{\pi_1^{temp}( X_{{\C}_p}^{^{an}},x)}&\to &\widehat{\pi_1^{temp}( X^{^{an}},x)}&\to &
G_K &&\\ && \downarrow \cong &&\downarrow \cong &&\downarrow  && \\ \{1\}&\to &
 {\pi_1^{alg}( X_{{\C}_p} ,x)}&\to & {\pi_1^{alg}( X,x)}&\to &
G_K &\to &\{1\} \end{matrix}\]
and the last row is the canonical exact sequence. In particular, the actions $\rho$ (\resp $\rho_x$) of 6.1.1
are compatible with the homonymous actions defined in 2.1.
 
\bigskip \section{  $Out\,\pi_1^{temp} \;$ for a curve}\label{s5}

\subsection{} In this section, we examine the group $Out\,
\pi_1^{temp}(X^{^{an}})$ when $X$ is a smooth algebraic curve over $\C_p$. 
 We fix a base point $x$. 

By {\it characteristic quotient} of a given topological group $G$, we mean the quotient of $G$ by some
characteristic open subgroup. 

\begin{lemme}  $\pi_1^{top}(X^{^{an}},x)$ is a characteristic
quotient of $\pi_1^{temp}(X^{^{an}},x)$.  
\end{lemme}

\noindent{\it Remark.} We do not know whether this holds for higher dimensional $X$, if
$\pi_1^{top}(X^{^{an}},x)$ has torsion.

\proof For a curve, $\pi_1^{top}(X^{^{an}},x)$ is a free, hence torsion-free, discrete group. By proposition
4.3.1, it can actually be described as the biggest torsion-free discrete quotient of
$\pi_1^{temp}(X^{^{an}},x)$. This group-theoretic characterization makes it clear that the kernel of
$\pi_1^{temp}(X^{^{an}},x)\to
\pi_1^{top}(X^{^{an}},x)$ is characteristic (and of course open).
\qed

 Let us recall from 4.5. that 
\[\pi_1^{temp}(X^{^{an}},x)= \varprojlim \,\Delta_n      \] 
 where $\Delta_n$ sits in an exact sequence
  \[\{1\}\to \pi_1^{top}(X_n^{^{an}})\to \Delta_n\to \Gamma_n\to \{1\}. \]  

\begin{lemme} $i)\,$ $\Delta_n$ is a characteristic quotient of  
$\pi_1^{temp}(X^{^{an}},x)$. 
\\ $ii)\,$ If $n\vert n'$, both $\Gamma_n$ and $\Delta_n$ are characteristic quotients of $\Delta_{n'}$. In
particular (for $n=n'$), the free group $\pi_1^{top}(X_n^{^{an}})$ is a characteristic subgroup of
$\Delta_n$.    
\end{lemme}

\proof $i)$ Let us consider again the commutative diagram from 4.5 \[\begin{matrix} \{1\}&\to
&\pi_1^{temp}(X_n^{^{an}})&\to &\pi_1^{temp}(X^{^{an}})&\to &\Gamma_n&\to &\{1\}&\\
\downarrow &&\downarrow &&\downarrow &&\downarrow &&\downarrow &\\ \{1\}&\to &\pi_1^{top}(X_n^{^{an}})&\to
&\Delta_n&\to &\Gamma_n&\to &\{1\}&. \end{matrix}\]
By definition, the kernel of $\widehat{\pi_1^{temp}(X^{^{an}})}= \pi_1^{alg}(X)\to\Gamma_n$ is
characteristic. This implies that the kernel $\pi_1^{temp}(X_n^{^{an}})$ of $ \pi_1^{temp}(X^{^{an}})
\to\Gamma_n$ is characteristic. On the other hand, the kernel of $\pi_1^{temp}(X_n^{^{an}} )\to
\pi_1^{top}(X_n^{^{an}})$ is characteristic, by 6.1.1. It follows that the kernel of $ \pi_1^{temp}(X^{^{an}})
\to \Delta_n$ is characteristic.

$ii)$: the profinite completion $\widehat{\Delta_{n'}}$ of $\Delta_{n'}$ is a quotient of $\pi_1^{alg}(X)$
which maps to $\Gamma_n$. By observation 2.2.1,
$\Gamma_n$ is a characteristic quotient of
$\widehat{\Delta_{n'}}$, hence of $\Delta_{n'}$; in other words, the kernel $K_{n'}$ of $\Delta_{n'}\to
\Gamma_n$ is characteristic.  

The kernel of
$K_{n'}\to
\pi_1^{top}(X_n^{^{an}})$ is the image in $\Delta_{n'}$ of the kernel of $ \pi_1^{temp}(X_{n}^{^{an}})
\to \pi_1^{top}(X_n^{^{an}})$. By the same argument as in 6.1.1 (using 4.3.1), it can be characterized as the
smallest normal subgroup $N$ of $K_{n'}$ such that $K_{n'}/N$ is torsion-free; this shows that it is
characteristic in $K_{n'}$, hence in $\Delta_{n'}$. It follows that  $ \Delta_n$ is a characteristic quotient
of  $\Delta_{n'}$.
  \qed 

\begin{cor} \[  Aut\,\pi_1^{temp}(X^{^{an}},x)= \varprojlim \,Aut \,\Delta_n   ,\;\; Out\,\pi_1^{temp}(X^{^{an}})=
\varprojlim \,Out \,\Delta_n .\]
 \end{cor}
Indeed, the previous lemma makes sense of the inverse limit $ \varprojlim \,Aut \,\Delta_n  $ and provides a
natural homomorphism $Aut\, \pi_1^{temp}(X^{^{an}},x)\to  \varprojlim \,Aut \,\Delta_n  $, which has an obvious
inverse (note however that the transition maps need not be surjective).\qed 

We endow $ Aut\,\pi_1^{temp}(X^{^{an}},x)$ and $ Out\,\pi_1^{temp}(X^{^{an}},x)$ with the inverse limit
topology. With this topology, it is clear that $ Aut\,\pi_1^{temp}(X^{^{an}},x)$ acts continuously on
$\pi_1^{temp}(X^{^{an}},x)=\varprojlim \,  \,\Delta_n $.

\begin{prop}  $Out\,
\pi_1^{temp}(X^{^{an}})$ is a
countable inverse limit of discrete finitely generated groups; in particular it is a polish group (\cf 4.4). 
 \end{prop}

\proof Indeed, finitely generated virtually free groups (such as $\Delta_n$) have finitely
generated\footnote{even finitely presented} automorphism group (Krstic, McCool \cite{mc}).\qed 
 
\subsection{} For a curve, we have seen that $\pi_1^{temp}(X^{^{an}},x)= \varprojlim  \,\Delta_n$ embeds into its
profinite completion $\pi_1^{alg}(X,x)=\varprojlim  \,\Gamma_n$. A fortiori, $Aut\,\pi_1^{temp}(X^{^{an}},x)=
\varprojlim\,  Aut\,\Delta_n$ embeds into $Aut\,\pi_1^{alg}(X,x)=\varprojlim  \,Aut\,\Gamma_n$. For outer
automorphisms, we also have:

\begin{prop} $Out\,\pi_1^{temp}(X^{^{an}} )$ embeds into $Out\,\pi_1^{alg}(X)\,$; in particular, it is
residually finite.
\end{prop}
\proof If $X$ is not hyperbolic, $\pi_1^{temp}(X^{^{an}})$ is abelian and coincides with $\pi_1^{alg}(X)$
except in the case of an elliptic curve with bad reduction, for which $\pi_1^{temp}(X^{^{an}})\cong \Z \times
\hat\Z$; in either case, the statement is trivial. 

\noindent We now assume that $X$ is hyperbolic. The assumption
of hyperbolicity ensures that the center of $\pi_1^{temp}(X^{^{an}},x)$ is trivial (an immediate consequence of
the fact that the center of its profinite completion $\pi_1^{alg}(X,x)$ is trivial).

\noindent Let us assume that  $\varpi\in \pi_1^{alg}(X,x)$ induces an automorphism of
$\pi_1^{temp}(X^{^{an}} ,x)$ by conjugation. We have to show $\varpi \in \pi_1^{temp}(X^{^{an}}
,x)$. 

\noindent Let us fix $n$, and consider the composed epimorphism
\[\pi_1^{alg}(X,x)\surj \widehat{\Delta_n} \surj \Gamma_n \]
(where $\widehat{\Delta_n}$ denotes the profinite completion of $\Delta_n$, which sits in an exact sequence
\[\{1\} \to \widehat{\pi_1^{top}(X_n^{^{an}} )}\to \widehat{\Delta_n} \to  \Gamma_n \to \{1\}\;\; ). \] Let $\pi_n$
be any element of
$\pi_1^{temp}(X^{^{an}} ,x)$ which has the same image as $\varpi$ in $\Gamma_n$. The image $\varpi_n$ of
$\varpi.\pi_n^{-1}$ in 
$\widehat{\Delta_n}$ actually belongs to $\widehat{\pi_1^{top}(X_n^{^{an}} )}$, and induces
by conjugation an automorphism of the free group $\pi_1^{top}(X_n^{^{an}} )$ (\cf 6.1.1, 6.1.2). By 3.2.1, it
follows that
$\varpi_n\in \pi_1^{top}(X_n^{^{an}}) $, hence the image of $\varpi$ in $\widehat{\Delta_n}$ actually belongs to
$\Delta_n$. Therefore $\varpi \in \pi_1^{temp}(X^{^{an}}
,x)= \varprojlim 
\,\Delta_n$.\qed
  
\begin{cor} If $X$ is hyperbolic, $\pi_1^{temp}(X^{^{an}},x)$ is its own normalizer in its profinite
completion. \qed
 \end{cor}
 
\noindent We denote by $\overline{ Out}\,{\pi_1^{temp}(X^{^{an}}  )}$ the closure of
$Out\,{\pi_1^{temp}(X^{^{an}}  )}$ in
$ Out\,\pi_1^{alg}(X)$: \[\overline{ Out}\,{\pi_1^{temp}(X^{^{an}}  )}=  
\varprojlim \,\Im\,[Out\,{\pi_1^{temp}(X^{^{an}}  )} \to Out\,\Gamma_n] \] 
(this limit is clearly a quotient of the profinite completion of $Out\,{\pi_1^{temp}(X^{^{an}}  )}$ and embeds
into $\varprojlim \, Out\,\Gamma_n =   Out\,\pi_1^{alg}(X)$ ). 

\medskip 
\noindent{\it Remark.} For the problems studied in the sequel (\cf 8.7.2), it would be interesting to know
whether
$Out\,{\pi_1^{temp}(X^{^{an}}  )}= \varprojlim 
\,Out\,\Delta_n $ is a profinite group or not, especially in the case $X={\bf
P}^1_{\C_p} \setminus\{0,1,\infty\}$. This involves the question of finiteness of the outer automorphism groups
of the finitely generated virtually free groups $\Delta_n$. In the abstract setting, a complete anwser to this
kind of question has been given by M. Pettet in terms of graphs of
groups \cite{pe}. His criteria apply for instance to Kato's presentation of
$p$-adic triangle groups as amalgams \cite{kat}, and could apply in the same way to the `$p$-adic
dessins d'enfants' invoked at the end of 4.5.  

Another interesting open problem is whether the natural homomorphism $Out\,{\pi_1^{temp}(X )}\to GL_2(\hat \Z)$
is surjective for $X={{\bf
P}^1_{\C_p}}^{^{an}}  \setminus\{0,1,\infty\}$ (one shows easily that the natural composed map
 \[{\pi_1^{temp}(X)}\to {\pi_1^{alg}(X)}\to H_1^{et}(X,\hat\Z)\cong\hat\Z^2\] is surjective).

\bigskip \section{  $G_\Q, G_{\Q_p}$ and geometric fundamental groups }\label{s6}

\subsection{}
 
Now $k$ is a number field, and we fix an embedding $\iota_v:\,\bar k\inj \C_p$. We denote by $k_v$ the
topological closure of $k$ in $\C_p$. 

For any smooth geometrically connected algebraic $k$-variety $X$, we denote by $X^{^{an}}$ the $p$-adic
analytic manifold attached to $X_{\C_p}$. 

Let us fix a $\C_p$-point $x$ of $X$. One has canonical isomorphisms 
\[\pi_1^{alg}(X_{\bar k},x) =
\pi_1^{alg}(X_{\C_p},x)=
\widehat{\pi_1^{temp}(X^{^{an}},x)},\]  the latter group being the profinite completion of the
tempered fundamental group ${\pi_1^{temp}(X^{^{an}},x)}$ which classifies tempered \'etale cove\-rings 
of $X^{^{an}}$.

There is an outer action $\rho$ of $G_k$ on $\pi_1^{alg}(X_{\bar k},x) $ and a compatible outer action of
$G_{k_v}$ on $\pi_1^{temp}(X^{^{an}},x)\;$ (6.1, 6.2).

\subsection{} Let us now assume that $X$ is a hyperbolic affine curve. We have
encountered two closed subgroups of
$Out\,\pi_1^{alg}(X_{\bar k})$: 
$\overline{Out}\,{\pi_1^{temp}(X^{^{an}}  )}$ and
$G_k$. As in 3.3, we shall study their intersection, starting with the case of ${\bf
P}^1 \setminus\{0,1,\infty\}$. The following is a non-archimedean analogue of 3.3.1.

\begin{thm} If $X= {\bf P}^1_k\setminus\{0,1,\infty\}\;$, then \[\;Out\,{\pi_1^{temp}(X^{^{an}}  )}\cap
\, G_k \,= \overline{ Out}\,{\pi_1^{temp}(X^{^{an}}  )}\cap
\, G_k  =G_{k_v}. \]
   \end{thm}
 
 \prf We already know that $G_{k_v}\subset  \;Out\,{\pi_1^{temp}(X^{^{an}} ,x )}\cap
\, G_k $ by 5.1.1 and 5.2. It remains to show that $\;\overline{ Out}\,{\pi_1^{temp}(X^{^{an}}  )}\cap
\, G_k \,\subset \, G_{k_v}  $. 

 We choose a rational base point $x$ of ${\bf P}^1_k\setminus\{0,1,\infty\}$.
 Let $\sigma\in G_k$ be such that $\rho(\sigma)\in \overline{ Out}\,{\pi_1^{temp}(X^{^{an}}  )}$. 
 We have to show that  
\[(\ast)\;\;\;\;\;{}^\forall j\in \bar k
\setminus
\{0,1\},\;\vert j\vert_p\leq 1\Rightarrow \vert j^\sigma\vert_p\leq 1.\]
 Let $E_j$ be an elliptic curve with invariant $j\in \bar k$. According to Belyi, there is a surjective
morphism
  $\varphi:E_j\to {\mathbf P}^1$ defined over $\Q(j)$ which is unramified outside
$\Sigma=\varphi^{-1}(\{0,1,\infty\})$.

 The idea of the proof is that if the image of some $\tau \in
Aut\,{\pi_1^{temp}(X^{^{an}} ,x )}$ in $ { Out}\,{\pi_1^{temp}(X^{^{an}}  )}$  is `sufficiently close' to
$\rho(\sigma)$, then $\tau$ induces an isomorphism $\pi_1^{temp}(E_{j}^{^{an}})\to
\pi_1^{temp}(E_{j^\sigma}^{^{an}}) ,$  whence $(\ast)$, since the algebraic structure of these fundamental
groups is determined by the alternative `good/bad' reduction, \cf 4.6.

\medskip Let us explain the details. We fix an integer $m$ so that the
characteristic finite etale covering $X_m$ of $X_{\bar k} $ defined in 2.2 dominates $
E_j\setminus
\Sigma $. This provides a sequence of finite etale coverings with base points
\[(X_m,x_m)\to (E_j\setminus
\Sigma,y) \to (X_{\bar k}, x) .\]
 One has a
commutative diagram
\[\begin{matrix}  \pi_1^{alg}(X_m,x_m) \;&  \stackrel{\rho_{x_m,x_m^\sigma}(\sigma)}\to &  \pi_1^{alg}(
X_m^\sigma,x_m^\sigma)\\
  \downarrow\;& &\downarrow   \\ \pi_1^{alg}(E_j\setminus
\Sigma,y) \;&  \stackrel{\rho_{y,y^\sigma}(\sigma)}\to &  \pi_1^{alg}( E_{j^\sigma}\setminus
{\Sigma}^\sigma, y^\sigma)\\
 \downarrow\;& &\downarrow   \\
               \pi_1^{alg}(X_{\bar k},x) \;& \stackrel{\rho_x(\sigma)}{\to} &  \pi_1^{alg}(X_{\bar
k},x).\\\end{matrix} \] 
 By assumption, \[\rho(\sigma)\in \overline{Out}\, \pi_1^{temp}(X^{^{an}} ) = \varprojlim \, \Im
[Out\, \pi_1^{temp}(X^{^{an}}) \to Out \,\Gamma_n]\]
 This implies the existence of $\tau\in Aut\,{\pi_1^{temp}(X^{^{an}} ,x )},
\;\varpi\in \,\pi_1^{alg}(X_{\bar k},x) $ such that  
$\rho_x(\sigma)$ and $\tau\circ ad(\varpi)$ have the same image in $\Gamma_m$. 
Since $\pi_1^{alg}(X_m,x_m)$ is open in the profinite completion $\pi_1^{alg}(X_{\bar k},x)$ of
${\pi_1^{temp}(X^{^{an}} ,x )}$, one has
\[\pi_1^{alg}(X_{\bar k},x)={\pi_1^{temp}(X^{^{an}} ,x )}.
\pi_1^{alg}( X_m,x_m),\] which allows to take $ \varpi\in \pi_1^{alg}( X_m,x_m)$.
Then there is a
similar commutative diagram built on $\tau$   
\[\begin{matrix}  \pi_1^{alg}(X_m,x_m) \;&  \to &  \pi_1^{alg}(
X_m^\sigma,x_m^\sigma)\\
  \downarrow\;& &\downarrow   \\ \pi_1^{alg}(E_j\setminus
\Sigma,y) \;&  \stackrel{\widehat{\tau'}}{\to} &  \pi_1^{alg}( E_{j^\sigma}\setminus
{\Sigma}^\sigma, y^\sigma)\\
 \downarrow\;& &\downarrow   \\
               \pi_1^{alg}(X_{\bar k},x) \;& \stackrel{\hat\tau}{\to} &  \pi_1^{alg}(X_{\bar
k},x),\\\end{matrix} \] 
which comes, by profinite completion, from a commutative diagram
\[\begin{matrix}  \pi_1^{temp}(X_m^{^{an}},x_m) \;&  \to &  \pi_1^{temp}(
{X^\sigma_m}^{^{an}},x_m^\sigma)\\
  \downarrow\;& &\downarrow   \\ \pi_1^{temp}(E_j^{^{an}}\setminus
\Sigma,y) \;&  \stackrel{\tau'}{\to} & \pi_1^{temp}( E_{j^\sigma}^{^{an}}\setminus
{\Sigma}^\sigma, y^\sigma)\\
 \downarrow\;& &\downarrow   \\
              \pi_1^{temp}(X^{^{an}},x) \;& \stackrel{\tau}{\to} & \pi_1^{temp}(X^{^{an}},x) \\\end{matrix} \] 
where horizontal maps are isomorphisms ( $\tau'$ stands for the induced isomorphism).

 On the other hand, the open immersion $E_j^{^{an}}\setminus \Sigma
\inj
E_j^{^{an}}$ induces a (strict) surjective homomorphism $ \pi_1^{temp}(E_j^{^{an}}\setminus
\Sigma,y) \to  \pi_1^{temp}(E_j^{^{an}}) $ \cite[\S\,4.5.5.b]{a}, which is nothing but the abelianization
homomorphism (this becomes clear on the profinite completions, or using 3.2.2.a). One thus has another
commutative diagram
\[\begin{matrix}   \pi_1^{temp}(E_j^{^{an}}\setminus
\Sigma,y) \;&  \stackrel{\tau'}{\to} & \pi_1^{temp}( E_{j^\sigma}^{^{an}}\setminus
{\Sigma}^\sigma, y^\sigma)\\
 \downarrow\;& &\downarrow   \\
               \pi_1^{temp}(E_j^{^{an}}
 ) \; & \stackrel{{\tau'}^{ab}}{\to}  & \pi_1^{temp}(E_{j^\sigma}^{^{an}} ) \\ \end{matrix} \] 
where horizontal maps are isomorphisms. 

 If
$\vert j\vert_p
\leq 1$,
$E_j$ has good reduction and it follows that $\pi_1^{temp}(E_j^{^{an}})=
\pi_1^{alg}(E_j)$, hence
$\pi_1^{temp}(E_{j^\sigma}^{^{an}})=
\pi_1^{alg}(E_{j^\sigma})$, which implies in turn that $E_{j^\sigma}$ has also good reduction at $p$, whence
$\vert j^\sigma\vert_p\leq 1$.
  \qed

\begin{cor} If $X= {\bf P}^1_k\setminus\{0,1,\infty\}\;$, the image of $\;Out\,{\pi_1^{temp}(X^{^{an}}  )}$ in 
$\;Out\,{\pi_1^{alg}(X_{\bar k} )} $ is not dense in any open subgroup.  \qed
\end{cor}

For more general affine hyperbolic $k$-curves $X$, we have the following non-archimedean analogue of 3.3.2:

\begin{thm} Assume that some etale cove\-ring of $X_{\bar k}$ admits a non-constant rational function which
omits at least three values.
  Then \[\;Out\,{\pi_1^{temp}(X^{^{an}}  )}\cap
\, G_k \,=     \, G_{k_v}.\] A fortiori, if the geometric
point $x$ is $k$-rational, then the subgroup of
$G_k$ which stabilizes
$\pi_1^{temp}(X^{^{an}} ,x )\subset \pi_1^{alg}(X_{\bar k}, x)\;$ is $\;G_{k_v}$.
\end{thm}

\prf From the assumption, there is a finite extension $k'$  of $k$ in $\C_p$,
a finite etale cove\-ring
$Y\to X_{k'}$ and a dominant morphism  \[Y\stackrel{\varphi}\to Z={\bf
P}^1_{k'}\setminus\{0,1 ,\infty\} .\] (We may also assume that $X$, $Y$, $Z$ are endowed with compatible
$k'$-rational geometric points $x, y,z$ respectively).
The result
will follow from 7.2.1 by the following d\'evissage:

\begin{lemme} The statement \[(\ast)_{X,k}\;\;\;\;\; { Out}\,{\pi_1^{temp}(X^{^{an}}  )}\cap \, G_k \,\subset \,
G_{k_v}  \] satisfies
\\$ a)$  $(\ast)_{X,k'}  \Rightarrow (\ast)_{X,k}\; $ if $k'$ is a finite extension of $k$ in $\C_p$  ,   
\\$ b)$ $(\ast)_{Z,k} \Rightarrow (\ast)_{Y,k} \;$ if there is a non-constant morphism $f:\,Y\to Z$.
\\$ c)$  $(\ast)_{Y,k}  \Rightarrow (\ast)_{X,k}\; $ if $Y$ is a finite \'etale cove\-ring of $X$,
\end{lemme}
\prf b) and c) are completely parallel to the proof of 3.3.2, replacing $\C$ by $\C_p$, $\pi_1^{top}$ by
$\pi_1^{temp}$, the reference to 3.2.1 by a reference to 6.2.2, and sublemma 3.3.4 by

\begin{slemme} (Neukirch, \cf
\cite[\S\,12.1.10]{nsw}) If $H$ is a closed subgroup of $G_k$ containing $G_{k_v}$ as a subgroup of finite index,
then $H=G_{k_v}$.\qed
\end{slemme}

 As for a):   
 let $H$ stand for $Out \,{\pi_1^{temp}(X^{^{an}}  )}\cap \, G_k $. It contains $G_{ k_v}$. On the other
hand,  $Out \,{\pi_1^{temp}(X^{^{an}}  )}\cap \, G_{k'} $, which is nothing but $G_{  {k'}_v}$ due to 
$(\ast)_{X,k'}$, is a subgroup of finite index of $H$. Therefore $H$ is closed in $G_k$ and contains
$G_{k_v}$ as a subgroup of finite index. One concludes by 7.2.4.\qed
  
\medskip\noindent {\it Remarks.} 1) One version of a celebrated theorem of Mochizuki \cite{mo'} says that for a
hyperbolic curve $X_{k_v}$ over $k_v$, the centralizer of $G_{k_v}$ in $Out\,\pi_1^{alg}(X_{\C_p}) $ is
isomorphic to $Aut\, X_{k_v}$. Denoting by $X^{^{an}}  $ the $p$-adic manifold attached to
$X_{\C_p}$, it follows that the center of $Out \, \pi_1^{temp}(X^{^{an}}  ) $ is isomorphic to the center of
$Aut\,X_{k_v}$. In particular, it is trivial if $X_{k_v}={\bf P}^1_{k_v}\setminus\{0,1,\infty\}$.
 
2) From $G_{k_v}$, one can recover the inertia group and the Frobenius elements by a purely group-theoretic
recipe
\cite[\S\, 12.1.8]{nsw}.
 
3) Theorems 7.2.1 and 7.2.3 do not change if one composes $\iota_v$ with any element of $G_{k_v}$. Thus
it is enough to fix a $p$-adic place $v$ of $\bar k$ (instead of $\iota_v$).

\bigskip \section{ The local Grothendieck-Teichm\"uller group at $p$}\label{s7}

\subsection{} In this section, we define a ``$p$-local'' version of the Grothendieck-Teichm\"uller group
$\widehat{GT}$, which plays with respect to $G_{\Q_p}$ the role that $\widehat{GT}$ plays with respect to
$G_\Q$.

\medskip  For  
$r\geq 4$, let 
$\sM_{0,r}$ be the moduli space for curves of genus
$0$ with
$r$ ordered marked points. This smooth geometrically connected $\Q$-variety carries an action of the symmetric
group
${\frak S}_r$ (permutation of the marked points). 
 It admits a canonical smooth compactification $\overline\sM_{0,r} $
(Deligne-Mumford-Knudsen), such that the boundary $\partial \sM_{0,r} $ is a divisor with normal crossings;
over $\C$, this compactification is simply-connected \cite[\S\, 1.1]{bp}, and a fortiori,
$\pi_1^{alg}(\overline\sM_{0,r, {\bar \Q}})=\{1\}$ (from now on, we drop the implicit base
point from the notation).

For $r=4$, $\sM_{0,4}\cong {\bf P}^1_\Q\setminus \{0,1,\infty\},\;\overline{\sM}_{0,4}\cong {\bf P}^1_\Q,\;$ and
the action of ${\frak S}_4$ factors through the standard action of ${\frak S}_3$.  
 In general, $\sM_{0,r}$ is isomorphic to the complement in $ ({\bf P}^1_\Q\setminus \{0,1,\infty\})^{r-3}$
of the partial diagonals defined by the equality of two coordinates.  

  Forgetting the last marked point provides an ${\frak S}_r$-equivariant surjective morphism
$p_{r}: \sM_{0,r+1} \to  \sM_{0,r} $ with fibers isomorphic to $ {\bf P}^1  $ minus $r$ points. 

\subsection{} The fundamental group
$\pi_1^{top}(\sM_{0,r,\C }^{^{an}})$ (a {\it mapping class group in genus zero}) is
generally denoted by
$\Gamma_{0,r}$, but we add a superscript
$\;{}^{top}\;$ in order to distinguish it from subsequent variants. 
 
The fibration $p_{r}$ gives rise to an exact sequence of
fundamental groups
\[ \{1\} \to F_{r-1} \to \Gamma^{top}_{0,r+1}\to \Gamma^{top}_{0,r}\to \{1\}\]
where $F_{r-1}\cong \pi_1^{top}({{\bf
P}^1_\C}^{^{an}}\setminus \{ x_1,\ldots, x_r\})$ is the discrete free group on $r-1$ generators. On the other
hand, according to a classical result of Nielsen-Zieschang (\cf \cite{zvc}), the corresponding homomorphism
\[\Gamma^{top}_{0,r}\to  
 Out\,F_{r-1}  \] induces an isomorphism
\[\Gamma^{top}_{0,r}\cong
Out^\natural\,\pi_1^{top}({{\bf P}^1_\C}^{^{an}}\setminus \{ x_1,\ldots, x_r\}) \cong \ker [Out\,F_{r-1} \to
GL_{r-1}(\Z)]\] (the sign $^\natural$ indicates the subgroup of (outer) automorphisms which fix the
conjugacy classes of the local monodromies at the points $x_i$). [From there, using the surjectivity
of $Aut\,F_n \to GL_n(\Z)$ \cite[\S\,I.4.4]{lys} , one derives the
structure of $Aut\,F_n$ (\resp
$Out\,F_n$) for any $n\geq 1$: it is an extension of $GL_n(\Z)$ by an iterated extension of free groups of
successive ranks $2,3,\ldots, n$ (\resp $2,3,\ldots, n-1$)].

 \subsection{} The {\it Grothendieck-Teichm\"uller group} $\widehat{GT}$ introduced by
Drinfel'd\footnote{unlike the usage in the present paper, $\widehat{GT}$ is just a
(traditional) notation and does not indicate the profinite completion of some group $GT$; this notation helps
distinguishing the profinite Grothendieck-Teichm\"uller group from other variants (notably the pro-unipotent
one related to mixed Tate motives)}
 in connection with the theory of quantum groups is a closed subgroup of the profinite group $Aut\,\hat F_2$. If
$x,y$ denote fixed topological generators of $\hat F_2$, it consists of automorphisms of the form $x\mapsto
x^\lambda,\, y\mapsto f^{-1}y^\lambda f$ for $\lambda \in \hat\Z^\ast$ and appropriate $ f\in \hat F_2'$
satisfying three well-known equations I, II, III, which we do not recall here (\cf \cite{dr} \cite{i}).

 The bigger group obtained by discarding equation III (the pentagonal equation) is sometimes denoted by
$\widehat{GT}_0$.  
  
  \subsection{} For any $r\geq 4$, we set $\Gamma^{alg}_{0,r}: = \pi_1^{alg}(\sM_{0,r ,\bar
\Q} )$. Thus $\Gamma^{alg}_{0,r}$ is identified with the profinite completion of $\Gamma^{top}_{0,r}$
if one fixes an embedding
$\iota : \bar\Q \inj \C$.  

Following Nakamura \cite{na}, we denote by
$Out^\flat\,\Gamma^{alg}_{0,r}\;$ the subgroup of outer automorphisms $\sigma$ of $ \Gamma^{alg}_{0,r}$ which
sends the conjugacy class of a generator of the inertia group at each component of $\partial
\sM_{0,r}
$ to the conjugacy class of the same generator raised to some power  
$ \lambda(\sigma) \in \hat\Z^\ast$ (the same $\lambda(\sigma)$ for each component). There is an obvious
homomorphism \[\lambda : Out^\flat\,\Gamma^{alg}_{0,r}\to \hat\Z^\ast, \;\sigma\mapsto \lambda(\sigma).\]
We denote by 
$\;\widehat{GT}^{(r)}:= Out^\flat_{\frak S_{r
}}\, \Gamma^{alg}_{0,r}\;$ the subgroup of elements which commute with
the action of ${\frak S_{r
}}$. This notation is motivated by: 
  
 \begin{thm}(Harbater, Schneps \cite{hs})
 
 \noindent 1) The ${\frak S}_r$-equivariant surjective morphism
$ \sM_{0,r+1}\stackrel{p_{r}}{\longrightarrow}   \sM_{0,r} $ induces canonical homomorphisms
 \[Out^\flat\,\Gamma^{alg}_{0,r+1}\to Out^\flat\,\Gamma^{alg}_{0,r},\;\;\;\; \widehat{GT}^{(r+1)}
\;\;\stackrel{q_{r}}{\rightarrow} \;\;\;\widehat{GT}^{(r)}  \]
 compatible with the $\lambda$-map.

 \noindent 2) For all $r\geq 5$, $q_{r}: \widehat{GT}^{(r+1)} \to \widehat{GT}^{(r )} $ is an
isomorphism, 

\noindent  3) An embedding $\iota : \bar\Q \inj \C$ being fixed, there is a natural commutative diagram
\[\begin{matrix} \widehat{GT}^{(5)} &\stackrel{q_{4}}{\inj} & \widehat{GT}^{(4)} \\ \downarrow \cong
&&\downarrow \cong\\ \widehat{GT}& \inj & \widehat{GT}_0  \end{matrix}\]  where the vertical maps are
isomorphisms.\qed 
 \end{thm} 
\noindent ($\iota$ is used to identify $\hat F_2$ with $ \Gamma^{alg}_{0,4}$ based at some tangential base
point).

  This result (partly based on results in \cite{na}) provides a `purely algebraic' description of $\widehat{GT}$
(independent of the topology of $\C$), and at the same time, a natural interpretation of the Drinfel'd-Ihara
embedding
$G_\Q\inj
\widehat{GT}$ (as an embedding of type $\rho$, with the notation of 2.1). The question of how close is
$\widehat{GT}$ to $G_\Q$ is a famous open problem.

\subsection{} In view of the results of \S 7, it is natural to try to replace the algebraic fundamental groups
which occur in the previous theorem by tempered fundamental groups in the $p$-adic context. 

We fix a prime number $p$, and consider the polish groups $\pi_1^{temp}(\sM_{0,r,\C_p }^{^{an}})$. It is likely
that these are residually finite, but having no proof, we are led to introduce the following $p$-adic avatar of
the mapping class group in genus zero:
 \[\Gamma^{temp}_{0,r}:=\pi_1^{temp}(\sM_{0,r,\C_p }^{^{an}})/\ker[\pi_1^{temp}(\sM_{0,r,\C_p }^{^{an}})\to
\widehat{\pi_1^{temp}(\sM_{0,r,\C_p }^{^{an}})}] .\]
 
\noindent Thus $\Gamma^{temp}_{0,r}$ is, like $\pi_1^{temp}(\sM_{0,r,\C_p }^{^{an}})$,  a polish group, and
embeds into its profinite completion, which can be identified with the
 $\Gamma^{alg}_{0,r}$ if one fixes an embedding $\iota_p: \bar\Q \inj \C_p$. 

\noindent Note that $\Gamma^{temp}_{0,4}$ is
our familiar group $ \pi_1^{temp}({{\bf P}^1_{\C_p}}^{^{an}}\setminus \{0,1,\infty\})$.

\begin{prop} The morphism $p_{r}:\,\sM_{0,r+1}\to \sM_{0,r}$ induces strict surjective homomorphisms
\[\pi_1^{temp}(\sM_{0,r+1,\C_p }^{^{an}})\surj \pi_1^{temp}(\sM_{0,r,\C_p }^{^{an}}),\;\;\;\;
\Gamma^{temp}_{0,r+1}\surj \Gamma^{temp}_{0,r}.\] 
\end{prop}

\prf From the description $\sM_{0,r,\C_p} = ({\bf P}^1_{\C_p}\setminus \{0,1,\infty\})^{r-3}\setminus \Delta
,$ we deduce a factorization of $p_{r}$ into a closed embedding $j:\,\sM_{0,r+1,\C_p}\inj {\bf A}^1_{\C_p}
\times\sM_{0,r,\C_p} $ followed by the second projection. This
second projection induces an isomorphism on $\pi_1^{temp}$: this follows from the fact that any finite etale
covering of ${\bf A}^1_{\C_p}
\times\sM_{0,r,\C_p} $ is of the form ${\bf A}^1_{\C_p}
\times\sM_1$ for some finite etale covering $\sM_1$ of $\sM_{0,r,\C_p} $  (in the algebraic or analytic
category, this does not matter thanks to \cite{lu}), and from the fact that any topological covering of 
${{\bf A}^1}_{\C_p}^{^{an}}\times\sM_1^{^{an}} $ is of the form ${{\bf A}^1}_{\C_p}^{^{an}}\times\sM_2^{^{an}}
$ for some topological covering $\sM_2^{^{an}}$ of
$\sM_1^{^{an}}$ (since $\pi_1^{top}({{\bf A}^1}_{\C_p}^{^{an}})=\{1\}$, and ${{\bf
A}^1}_{\C_p}^{^{an}}\times\sM_1^{^{an}} $ is homotopy equivalent to the naive cartesian product of 
${{\bf A}^1}_{\C_p}^{^{an}}$ and
$\sM_1^{^{an}}$ \cite[prf. of 10.1]{be''}).

 On the
other hand,
$j$ is the inclusion of the complement of $r$ smooth hyperplanes, and it follows from \cite[\S \,4.5.5.b)]{a}
that it induces a strict surjection on
$\pi_1^{temp}$ (this result uses heavily the fact that one is dealing with polish groups). 
\qed

\noindent {\it Remark.} One has $\pi_1^{temp}(\overline{\sM}_{0,r ,\C_p }^{^{an}})=\{1\}$: indeed, we have
seen that any finite etale covering of $\overline{\sM}_{0,r ,\C_p }$ is trivial (in the algebraic or analytic
category \cite{lu}), and on the other hand, since $\overline{\sM}_{0,r}$ is 
Zariski open in the affine space ${\bf A}^{r-3}$, $\pi_1^{top}(\overline{\sM}_{0,r ,\C_p
}^{^{an}})=\pi_1^{top}({{\bf A}_{\C_p}^{r-3}} )=\{1\}$ \cite[\S \, 1.1.4]{a}.

This suggests that $\pi_1^{temp}({\sM}_{0,r ,\C_p }^{^{an}})$ might `come' in some sense from neighborhoods of
points of maximal degeneration in $\partial {\sM}_{0,r ,\C_p }^{^{an}}$. 

\subsection{} For any $r\geq 4$, we denote by $Out^\flat\,\Gamma^{temp}_{0,r}\;$ the subgroup of outer
automorphisms $\sigma$ of
$ \Gamma^{temp}_{0,r}$ which send the conjugacy class of the local monodromy at each component of
$\partial \sM_{0,r} $ to the conjugacy class of the same local monodromy raised to some power  
$\lambda(\sigma) \in \hat\Z^\ast$ (the same $\lambda(\sigma)$ for each component); we refer to
\cite[\S\,4.5.4]{a} for the precise definition of this conjugacy class.   

\begin{notn} $\bullet$ We denote by $GT^{(r)}_p:= Out^\flat_{\frak S_{r }}\,
\Gamma^{temp}_{0,r}$ the subgroup of $Out^\flat\,
\Gamma^{temp}_{0,r}$ of elements which commute with the action of the symmetric group $\frak S_{r }$. 

\noindent $\bullet$ We denote by
$\widehat{GT}^{(r)}_p$ the closure of the image of
$GT^{(r)}_p$ in
$\widehat{GT}^{(r)}\cong 
\widehat{GT}$.   
 \end{notn}
\noindent (This notation is consistent with the notation $\widehat{GT}$ of
the Grothendiek-Teich\-m\"ul\-ler group, but it should not be confused with the profinite completion
$\widehat{GT^{(r)}_p}$ of
$GT^{(r)}_p$; we do not know whether the natural epimorphism $\widehat{GT^{(r)}_p} \to
\widehat{GT}^{(r)}_p$ is an isomorphism, but \cf 8.7.2.). 

Note that $\widehat{GT}^{(4)}_p$ is a closed subgroup of $\overline{Out}\,\pi_1^{temp}({{\bf
P}^1_{\C_p}}^{an}\setminus \{0,1,\infty\})$.

\begin{prop}  \noindent 1) The ${\frak S}_r$-equivariant surjective morphism
$  p_{r}  $ induces canonical {\rm injective} homomorphisms 
 \[  {GT}^{(r+1)}_p
\;\;\stackrel{q_{r}}{\inj} \;\;\; {GT}^{(r)}_p,\;\;\;\; \widehat{GT}^{(r+1)}_p
\;\;\stackrel{q_{r}}{\inj} \;\;\;\widehat{GT}^{(r)}_p  .\]
 
 \noindent 2) The canonical homomorphism $ {GT}^{(r)}_p \to \widehat{GT}^{(r)}_p $ is injective; in
particular, $ {GT}^{(r)}_p$ is residually finite.  \end{prop}

\prf The part of 1) concerning $\widehat{GT}^{(r)}_p $ follows directly from 8.4.1.
 From the embedding
$\Gamma^{temp}_{0,r} \inj \Gamma^{alg}_{0,r}$, we derive an embedding
  $Aut^\flat\,\Gamma^{temp}_{0,r}\inj Aut^\flat\,\Gamma^{alg}_{0,r}$. The ``Aut" version of
8.4.1 (\cf \cite[3.1]{hs}) says that $p_{r}$ induces an injective
homomorphism $Aut^\flat\,\Gamma^{alg}_{0,r+1}\inj Aut^\flat\,\Gamma^{alg}_{0,r}$.
By application of the sorite 3.2.2.a) to $G=\Gamma^{temp}_{0,r+1},\;H=\Gamma^{temp}_{0,r}$ and $
\phi=p_{r,\ast}$, we derive that $p_{r}$ induces an injective homomorphism $Aut^\flat\,\Gamma^{temp}_{0,r+1}\inj
Aut^\flat\,\Gamma^{temp}_{0,r}$. Finally, using the surjectivity of $\Gamma^{temp}_{0,r+1} \surj
\Gamma^{temp}_{0,r}$ (8.5.1), we see that $p_{r}$ induces in turn injective homomorphisms
\[Out^\flat\,\Gamma^{temp}_{0,r+1}\inj Out^\flat\,\Gamma^{temp}_{0,r},\;\;{GT}^{(r+1)}_p
  {\inj}   {GT}^{(r)}_p,\] which proves 1).

For 2), we compose the maps $q_{r}$ up to $r=4$ and get an injective homomorphism ${GT}^{(r)}_p
  {\inj}   {GT}^{(4)}_p \subset Out\, \Gamma^{temp}_{0,4} $. But we know (6.2.1) that $Out\, \Gamma^{temp}_{0,4}
$ embeds into $Out\, \Gamma^{alg}_{0,4} $. Since the composed map ${GT}^{(r)}_p
  {\inj} Out\, \Gamma^{alg}_{0,4} $ factors through $\widehat{GT}^{(r)}_p$, assertion 2) follows. 
\qed 

 \begin{defn} $\bullet$ We define $GT_p$ to be the intersection of the
images of the ${GT}^{(r)}_p$ in ${GT}^{(4)}_p$ (for all $\;r\geq 4$).

 \noindent $\bullet$ We define the {\it local Grothendieck-Teichm\"uller group at $p$} to be the closure
$\widehat{GT}_p$ of $GT_p$ in $\widehat{GT}^{(4)}_p$.
 \end{defn}

\noindent  From 8.6.2.1), we get 
\begin{cor} The abstract group $GT_p\;$ acts on the `tower' of tempered mapping class groups $
(\Gamma^{temp}_{0,r})_{r\geq 4}\,$.
\qed\end{cor} 

 On the other hand, $\widehat{GT}_p$ can be viewed as a
closed subgroup of the Grothendieck-Teichm\"uller group $\widehat{GT}$ via the Harbater-Schneps isomorphism $
\widehat{GT}^{(4)}\cong 
\widehat{GT}_0$.

\subsection{} We fix a $p$-adic place $v$ of $\bar \Q$. Then the maps $\rho$ from \S 2 and \S 5 give rise to
commutative squares of injective homomorphisms 
\[\begin{matrix} G_{\Q_p} &\to &{GT}_p &\;\;\;\;&&&& G_{\Q_p} &\to &{GT}^{(r)}_p\\ \downarrow &&\downarrow
&\;\;\;&&&& 
\downarrow &&\downarrow\\G_{\Q } &\to &\widehat{GT}_0 &\;\;\;\;&&&& G_{\Q} &\to &\widehat{GT}^{(r)} 
\end{matrix}\] which are compatible via the embeddings
$q_r$. It turns out that these squares are cartesian:

\begin{thm} $G_{\Q_p}=G_\Q \cap  {GT}_{p}= G_\Q \cap
\widehat{GT}_{p}  \; $ in $\;\widehat{GT}_0 \,$. Moreover,
for any $r\geq 4$, $G_{\Q_p}=G_\Q \cap  {GT}^{(r)}_{p} =G_\Q \cap \widehat{GT}^{(r)}_{p} \; $ in
$\;\widehat{GT}^{(r)} \;$.
\end{thm}
\proof This follows immediately from 7.2.1, since \[\widehat{GT}^{(4)}_p \subset
\overline{Out}\,\pi_1^{temp}({{\bf P}^1_{\C_p}}^{an}\setminus \{0,1,\infty\}). \;\;\;\;\;\;\hbox{ \qed} \] 

In the following corollary, we identify all $\widehat{GT}^{(r)}$ with $\widehat{GT}$ for $r\geq 5$ (according
to 8.4.1).

\begin{cor} If $G_\Q=\widehat{GT}$, then all these groups $GT_p= 
{GT}^{(r)}_{p}= 
\widehat{GT}_{p}=  \widehat{GT}^{(r)}_{p} \; $ would coincide with $G_{\Q_p}$ in
$\;\widehat{GT} $, for $r\geq 5$. 

In particular
$\;{GT}^{(5)}_{p}\, = \, Out^\flat_{\frak S_{5 }}\,
\Gamma^{temp}_{0,5}\;$  would be a finitely generated profinite group\footnote{topologically generated
by three elements, \cf
\cite[\S\,7.4.1]{nsw})}.\qed
\end{cor}

It is tempting to try to put this property in the wrong (\cf the remark after 6.2.2)...

\medskip
\noindent {\it Remarks.} 1) The analogue for $\widehat{GT}_p$ of the inertia subgroup of $G_{\Q_p}$ is the
kernel of the composed map  $\widehat{GT}_p \stackrel{\lambda}{\to} \hat\Z^\ast= ({\hat \Z}^{(p')})^\ast \times
\Z_p^\ast\to ({\hat \Z}^{(p')})^\ast .$ 

2) The image of $G_{\Q_p}$ by $\lambda$ is
$p^{\hat\Z}\times  \Z_p^\ast $. A natural question is whether the same holds for the image of $\widehat{GT}_p$
or even $\widehat{GT}^{(4)}_p$. After submission of this paper, A. Tamagawa has indicated to the author an
argument according to which this should follow from the main results of \cite{t}, notably theorem 5.2 and
proposition 5.3 of {\it loc. cit}.

3) Since $\lambda$ takes its values in $\hat{\Z}^\ast$, $GT^{(r)}_p$ acts on the quotient of 
$\Gamma^{temp}_{0,r}$ by the closed normal subgroup generated by some fixed integral powers of the local
monodromies at the components of the boundary. This quotient is an orbifold fundamental group in the sense of
\cite{a}, and its discrete representations correspond to certain algebraic vector bundles on $\sM_{0,r,\C_p}$
with fuchsian integrable connection and prescribed exponents \cite[\S 4.5.9]{a}, which form a $p$-adic analytic
subvariety of the corresponding algebraic moduli space of connections. In particular, $G_{\Q_p}$ acts through a
finite group on each of these $p$-adic analytic subvarieties.
 
4) As in the remark at the end of \S 7, one derives from Mochizuki's theorem (and 8.6.2) that the center of
any of the groups $GT_p, 
{GT}^{(r)}_{p}, 
\widehat{GT}_{p},  \widehat{GT}^{(r)}_{p} \; $ is trivial.

5) In complete analogy with the above, one can introduce the archimedean avatars $ {GT}^{(r)}_{\infty}$,
$\widehat{GT}^{(r)}_{\infty}$, ${GT}_{\infty}$, $\widehat{GT}_{\infty}$ (replacing $\C_p$ by $\C$ and
$\pi_1^{temp}$ by $\pi_1^{top}$). 
It turns out that all {\it these groups coincide with $G_\R$ for any $r\ge 4$}.  This is proved by the same
argument as in part $i)$ of the proof of 3.3.1, replacing the final reference to \cite[\S\,12.1.4]{nsw} by a
reference to
 \cite[prop. 4.ii]{ls}: $G_\R$ is its own centralizer in $\widehat{GT}$. 

6) Drinfel'd's presentation of $\widehat{GT}$ uses topological generators $x, y$ of $\hat F_2$ which are
the standard local monodromies at $0$ and $1$ in $\pi_1^{top}({{\bf P}^1_{\C}}^{an}\setminus
\{0,1,\infty\})$ (based at a suitable tangential base point). One could ask whether there is an analogous 
presentation of $\widehat{GT}$ of $p$-adic flavour, using local
monodromies at $0$ and $1$ in $ \pi_1^{temp}({{\bf P}^1_{\C_p}}^{an}\setminus
\{0,1,\infty\})$. It turns out that the answer is negative, as a consequence of the following observation: using
techniques of graphs of groups, we have shown in \cite[\S\,6.4.6]{a} that for $p=2$, one cannot find local
monodromies
$x,y,z$ at $0,1,\infty$ respectively such that $xyz = 1$.

 \section{ $\; G_{\Q_p} $ is the automorphism group of the tempered $\pi_1$-functor}\label{s8}

\subsection{} In this section, we deal with automorphisms of fundamental groups considered as functors.

\noindent Let $\sV_k$ be the category of all smooth geometrically connected algebraic varieties
defined over
$k$.  

\noindent Let $\sT $ be the category of separated topological groups up to inner automorphism: a morphism in
this category is a continuous homomorphism $G\to H$ given up to conjugation by an element of $H$ (in
particular, $Aut_\sT (G)= Out\,G$).

 In this article, we have been mainly concerned with the following functors \[\underline {\pi}_{\bar
k}^{alg}, \;\; \underline{\pi}_{\C_p}^{temp}\;: \;\;\sV_k\to \sT\]
\[\underline{\pi}_{ \bar k}^{alg}(X):= {\pi}_1^{alg}(X_{\bar k},x), \;\;\underline
{\pi}_{\C_p}^{temp}(X):=  {\pi}_1^{temp}(X_{\C_p}^{^{an}},x),\]

\noindent given an embedding $k\inj \C_p$ (note that in the category $\sT$, these fundamental groups do not
depend on the choice of $x$ up to {\it canonical} isomorphism; hence the above functors are uniquely defined up
to unique isomorphism; similarly, $\underline{\pi}_{\C_p}^{temp}$ depends only on the
$p$-adic place $v$ of
$k$ induced by the embedding, up to unique
isomorphism).

Let $k_v$ be the topological closure of $k$ in $\C_p$. Functoriality of the maps
$\rho$ considered in
\S 2 and
\S 5 may be expressed as follows: there are canonical homomorphisms 
\[ G_{k} \stackrel{\rho}{\to }  Aut\, \underline {\pi}_{  \bar k}^{alg}
 ,\;\;\; G_{k_v}  \stackrel{\rho}{\to }  Aut\,
\underline {\pi}_{\C_p}^{temp}   \] which build a commutative square if $\bar k$ is {\it the} closure
of $k$ in $\C_p$:
\[\begin{matrix} G_{k_v} &\stackrel{\rho}{\to } &Aut\,
\underline {\pi}_{\C_p}^{temp} & \\ \downarrow &&\downarrow
&\;\;\; \\G_{k} &\stackrel{\rho}{\to } &Aut\, \underline {\pi}_{  \bar k}^{alg} &. \end{matrix}\]
 
\subsection{} The following result had been conjectured by Oda and Matsumoto \cite{m'}:

\begin{thm}(Pop \cite{p}, unpublished) $\; G_{k}  \stackrel{\rho}{\to }  Aut\, \underline {\pi}_{ 
\bar k}^{alg} $ is an isomorphism.\qed
\end{thm}

\noindent (In fact, Pop proves more: $G_k$ coincides with  $Aut\, {\underline {\pi}_{ 
\bar k}^{alg}\,}_{\mid{\sW}}$ for any big enough full subcategory $\sW$ of $\sV_k$; we do not know any minimal
choice for $\sW$, but it suffices that $\sW$ contains the open subsets of the projective spaces). 

\medskip As a consequence of Pop's theorem and theorem 7.2.1, we get

\begin{thm} $\; G_{k_v}  \stackrel{\rho}{\to }  Aut\, \underline {\pi}_{\C_p}^{temp}$ is an
isomorphism. 
\end{thm}

\prf One has a sequence of homomorphisms
\[G_{k_v}  \stackrel{\rho}{\to }  Aut\, \underline {\pi}_{\C_p}^{temp}\,  {\begin{matrix} & Aut\, \underline
{\pi}_{ 
\bar k}^{alg} \stackrel{9.2.1}{=} G_k & \\ \nearrow &&\searrow \\ \searrow &&\nearrow \\  &
{Out}\,\pi_1^{temp}({{\bf P}^1_{\C_p}}^{an}\setminus \{0,1,\infty\}) \end{matrix}} 
 \,  {Out}\,\pi_1^{alg}({\bf P}^1_{\bar k} \setminus \{0,1,\infty\}) 
\] where the two maps on the right are injective. By 7.2.1, the intersection of $G_k $ and 
${Out}\,\pi_1^{temp}({{\bf P}^1_{\C_p}}^{an}\setminus \{0,1,\infty\}) $ in $ {Out}\,\pi_1^{alg}({\bf P}^1_{\bar
k} \setminus
\{0,1,\infty\}) $ is $G_{k_v}$. It follows that $\rho:\; G_{k_v}  {\to }  Aut\, \underline
{\pi}_{\C_p}^{temp}$ admits a left inverse $\nu$. It remains to show that $\,\ker \nu $ is trivial.

 Note that
$\ker \nu =\ker [Aut\, \underline{\pi}_{
\C_p}^{temp}\to Aut\, \underline {\pi}_{ 
\bar k}^{alg} ]$. It follows, by 6.2.1, that the evaluation of any element of $\,\ker \nu $ on
any curve
$X$ is the identity (as outer automorphism of the fundamental group). In order to reduce the case of an
arbitrary smooth geometrically connected $k$-variety $V$ to the case of a curve $X$, we choose a
quasi-projective dense open subset $U\subset V$ and an embedding of $U$ in some projective space ${\bf
P}^N_k$. According to a theorem of Bertini-Deligne \cite[1.4]{de}, there is a dense open subset $\sU$ of the
grassmannian of linear varieties of codimension $\dim U-1$ in ${\bf P}^N_k$ such that for all $L\in \sU(k)$,
the curve $X=L\cap U$ is geometrically connected and smooth, and the homomorphism 
\[\pi_1^{top}(X_\C^{^{an}})\to \pi_1^{top}(U_\C^{^{an}}) \] is surjective (for any fixed complex embedding
$k\inj
\C$). It follows that \[\pi_1^{alg}(X_{\bar k} )\to \pi_1^{alg}(U_{\bar k} ) \] is also
surjective; on the other hand
\[\pi_1^{alg}(U_{\bar k})\to \pi_1^{alg}(V_{\bar k}) \] is also surjective. Since the evaluation of
any element of $\,\ker \nu $ on
$X$ is the identity, the evaluation of
any element of $\,\ker \nu $ on
$V$ is also the identity, hence $\ker \nu=\{1\}$.
 \qed

 \bigskip 
 
\bigskip
 {\it Acknowledgements.} 
{A suggestion by Leila Schneps is at the origin of this paper: when I asked her, a few years ago,
about potential uses of tempered fundamental groups in Grothendieck-Teichm\"uller theory, she
immediately pointed out the problem of characterizing geometrically the local Galois groups. I am very
grateful for this pertinent suggestion; I also thank her, as well as Pierre Lochak, for many useful discussions
about Grothendieck-Teichm\"uller theory. Finally, I am much indebted to A. Tama\-gawa for his comments on a
first version of this text.}

\end{document}